\documentclass[11pt,a4paper]{amsart}
\usepackage{amsmath,amssymb,amsthm,mathtools}
\usepackage{booktabs}
\usepackage{enumitem}
\usepackage[margin=2.7cm]{geometry}
\usepackage{xcolor}
\definecolor{deepteal}{RGB}{0,128,128}
\usepackage[colorlinks=true,linkcolor=deepteal,citecolor=blue,urlcolor=deepteal]{hyperref}
\newif\ifshowchanges
\showchangesfalse
\newcommand{\chg}[1]{\ifshowchanges\textcolor{blue}{#1}\else#1\fi}

\numberwithin{equation}{section}
\newtheorem{theorem}{Theorem}[section]
\newtheorem{proposition}[theorem]{Proposition}
\newtheorem{lemma}[theorem]{Lemma}
\newtheorem{corollary}[theorem]{Corollary}
\theoremstyle{definition}
\newtheorem{definition}[theorem]{Definition}

\theoremstyle{remark}
\newtheorem{remark}[theorem]{Remark}
\newtheorem*{remark*}{Remark}

\newcommand{\Fq}{\mathbb{F}_q}
\newcommand{\A}{\mathbb{F}_q[T]}
\newcommand{\Kinf}{K_\infty}
\newcommand{\M}{\mathcal{M}}
\newcommand{\PP}{\mathcal{P}}
\newcommand{\V}{\mathcal{V}}
\newcommand{\F}{\mathcal{F}}
\newcommand{\G}{\mathcal{G}}
\newcommand{\I}{\mathcal{I}}
\newcommand{\nrm}[1]{\lVert #1\rVert}
\newcommand{\Wp}{W'_L}

\title[The distribution of $\alpha\pi$ modulo one in $\Fq(T)$]{The distribution of $\alpha\pi$ modulo one in $\Fq(T)$: beyond Matom\"aki's exponent}
\author{Esrafil Ali Molla}
\address{Department of Mathematics, Indian Institute of Technology Bombay, Mumbai 400076, India}
\email{esrafil.math@gmail.com}
\date{\today}

\begin{document}

\begin{abstract}
Let $\alpha$ be an element of $\Fq((1/T))$ which is not in $\Fq(T)$. We ask how small $\nrm{\alpha\pi}$ can be made when $\pi$ runs over irreducible polynomials, and write $\theta(\alpha)$ for the largest exponent $\theta$ such that $\nrm{\alpha\pi}<|\pi|^{-\theta}$ has infinitely many solutions. The analogue of Matom\"aki's exponent~\cite{Matomaki}, $\theta(\alpha)\ge\frac13$, was proved by Baier and the author~\cite{BM}, and it is the natural limit of \chg{methods based on} the Riemann Hypothesis. In this paper we go beyond it. For odd characteristic $p$ and $q$ large in terms of $p$ and $\theta$, we prove that $\theta(\alpha)\ge\theta$ for every $\alpha$ and every $\theta<\sqrt2-1=0.414\ldots$, and for every $\theta<\frac12$ when $\alpha$ has bounded partial quotients. The main tools are an exact count of points in affine subspaces of $\Fq[T]$, a splitting of these subspaces into long arithmetic progressions, and the results of Sawin and Shusterman~\cite{SS,SS2} on the M\"obius function along such progressions.
\end{abstract}

\maketitle

\tableofcontents

\section{Introduction}\label{sec:intro}

Throughout this article, $q$ is a power of a prime $p$, and $\varepsilon$ and $\eta$ denote small positive real numbers, fixed in each statement.

A fundamental theorem in Diophantine approximation is the following theorem due to Dirichlet.

\begin{theorem}[Dirichlet]\label{thm:dirichlet}
Given any real irrational $\alpha$, there are infinitely many pairs $(a,m)\in\mathbb Z\times\mathbb N$ of relatively prime integers such that
\[
\Bigl|\alpha-\frac am\Bigr|<m^{-2}.
\]
\end{theorem}

Equivalently,
\begin{equation}\label{eq:dirichlet}
\nrm{m\alpha}<m^{-1}
\end{equation}
for infinitely many positive integers $m$, where $\nrm{x}$ denotes the distance from $x$ to the nearest integer. This result is easy to prove using the pigeonhole principle or the continued fraction expansion of $\alpha$; the latter gives explicit solutions of the above inequalities.

It is natural to study the solubility of~\eqref{eq:dirichlet} when $m$ is restricted to subsets of the positive integers. The set of primes is a particularly interesting subset. The question now becomes: for which $\theta>0$ can one establish the existence of infinitely many primes $p$ such that
\begin{equation}\label{eq:primeineq}
\nrm{p\alpha}<p^{-\theta}\ ?
\end{equation}
In other words, we are interested in good approximations of $\alpha$ by fractions with prime denominator. This problem has a long history, and it has stimulated the development of several important tools in analytic number theory. We give a brief overview of this history.

In the thirties of the last century, I.~M.~Vinogradov~\cite{Vinogradov} established the existence of infinitely many primes $p$ satisfying~\eqref{eq:primeineq} with $\theta=\frac15$, by a sophisticated non-trivial treatment of trigonometric sums over primes, which also allowed him to prove the ternary Goldbach conjecture for sufficiently large odd integers without assuming any hypothesis like the Riemann Hypothesis. In 1977, Vaughan~\cite{Vaughan} improved the exponent $\frac15$ to $\frac14$, using his famous identity for the von Mangoldt function and refined Fourier analytic arguments. This was improved by Harman~\cite{Harman83} to $\frac3{10}$ in 1983, using his sieve method, which has since become an important tool in sieve theory~\cite{Harman}. There were several further improvements by Jia and Harman, who refined the numerical calculations and the sieve method itself; in particular, Harman~\cite{Harman96} obtained the exponent $\frac7{22}$. New ideas enabled Heath-Brown and Jia~\cite{HBJ} to come very close to the exponent $\frac13$ (they established $\theta=\frac{16}{49}$), and $\theta=\frac13$ was eventually achieved by Matom\"aki~\cite{Matomaki}. In particular, Heath-Brown and Jia made a connection to Kloosterman sums, and Matom\"aki utilised estimates for sums of Kloosterman sums, as established by Deshouillers and Iwaniec in their celebrated work~\cite{DI}. The exponent $\frac13$ is considered to be the limit of the current technology, and \chg{it is not hard} to establish it under the Generalised Riemann Hypothesis (GRH) for Dirichlet $L$-functions (see the comments in~\cite{HBJ}, and~\cite[\S2]{BM} for a short proof).

The problem has also been studied over number fields. Harman~\cite{HarmanG} proved the analogue of his exponent $\frac7{22}$ for the Gaussian primes, Baier~\cite{B2017} studied Diophantine approximation on lines in $\mathbb C^2$ with Gaussian prime constraints, and Baier and Mazumder~\cite{BMaz} treated real quadratic fields of class number one, where the unit group is infinite and the problem becomes two-dimensional. Assuming GRH for Hecke $L$-functions, Baier, Das and the author~\cite{BDM} extended Matom\"aki's exponent $\frac13$ to all quadratic number fields.

In function fields, GRH is known to hold: it is a theorem of Weil. Baier and the author~\cite{BM} translated the GRH proof of Matom\"aki's result into a version for the function field $\Fq(T)$ and obtained the analogue of the exponent $\frac13$ unconditionally, for odd prime powers $q\ge7$; for $q=4$ and $q=5$ they obtained the weaker exponent $1-\log q/\log\frac{q(q-1)}2$ (that is, $0.226\ldots$ and $0.301\ldots$), and for $q=2,3$ no non-trivial result. They also worked out a generalisation to imaginary quadratic extensions of $\Fq(T)$. Real quadratic extensions turned out to be harder, for the same reason as over number fields: the unit group is infinite. Baier and the author~\cite{BM2} first obtained the analogue of Vaughan's exponent $\frac14$ for real quadratic function fields of class number one, using a function field version of Vaughan's identity. More recently, the author~\cite{Mo} extended the exponent $\frac13$ to all real quadratic extensions of $\Fq(T)$, without any restriction on the class number, by writing the condition on normalised generators of principal prime ideals exactly as a finite combination of Hecke characters. Thus, over $\Fq(T)$ and its quadratic extensions, the exponent $\frac13$ is now available in complete generality. (The study of irreducible polynomials in residue classes, which underlies all these arguments, goes back to Hayes~\cite{Hayes}.)

In all these arguments the Diophantine condition is converted into a condition on residue classes modulo a convergent denominator $f$ of $\alpha$, and the Riemann Hypothesis gives square-root cancellation for each character modulo $f$ separately. This handles moduli up to $|f|\le X^{2/3}$, where $X$ is the size of the primes being counted. To cross this barrier, one needs information about primes, or about the M\"obius function, that cannot be extracted from individual $L$-functions.

In this article we go beyond the exponent $\frac13$ over $\Fq(T)$.\footnote{By ``Matom\"aki's exponent'' in the title we mean its function field analogue, the bound $\theta(\alpha)\ge\frac13$ over $\Fq(T)$ discussed above, rather than the exponent $\frac13$ of~\cite{Matomaki} over $\mathbb Z$.} The new input is a body of results on the M\"obius function in characteristic $p$ due to Sawin and Shusterman~\cite{SS,SS2}, which has no counterpart over $\mathbb Z$. Roughly speaking, they observed that in odd characteristic $\mu(r+s^p)$ behaves like a quadratic Dirichlet character in $s$, and this allowed them to prove cancellation of $\mu$ along arithmetic progressions far beyond the reach of the Riemann Hypothesis, provided $q$ is large compared with $p$. Our results may be summarised as follows, in terms of the exponent $\theta(\alpha)$ defined in \S\ref{sec:ffsetting}; the precise statements are given in \S\ref{sec:mainresults}.
\begin{itemize}[leftmargin=2em]
\item For odd $p$, every $\theta<\sqrt2-1=0.414\ldots$ and $q$ large in terms of $p$ and $\theta$, \emph{every} $\alpha$ satisfies $\theta(\alpha)\ge\theta$ (Corollary~\ref{cor:all}). In the language of the classical problem, this is the exponent $\frac13+\chg{\delta}$ for every $\chg{\delta}<\sqrt2-\frac43=0.080\ldots$; for instance, $\frac13+\frac1{126}$ is reached as soon as $q>13026\,p^2$ (Table~\ref{tab:q}).
\item If $\alpha$ has bounded partial quotients, or more generally if it has infinitely many scales at which two consecutive partial quotients are small, then $\theta(\alpha)$ is at least $\frac12-\eta$ (Corollaries~\ref{cor:tau} and~\ref{cor:sigma}).
\item If infinitely many convergent denominators of $\alpha$ are squarefree, then a theorem of Sawin~\cite{Sawin} gives $\theta(\alpha)\ge2\lambda(q)-1$, where $\lambda(q)\to1$ as $q\to\infty$, and this holds for every prime power $q$, including prime $q$ (Theorem~\ref{thm:sqfree}).
\end{itemize}
We note that the proofs of our main results (Theorems~\ref{thm:main} and~\ref{thm:main2}) do \emph{not} split into cases according to whether the moduli are squarefree; the results of Sawin and Shusterman that we use hold for arbitrary moduli. The squarefree case is treated separately only because, when it is available, a much stronger result follows by a different route, without any condition relating $q$ to $p$. 

\subsection{The function field setting}\label{sec:ffsetting}

Throughout, $q$ is a power of a prime $p$, $\A$ is the polynomial ring over the finite field $\Fq$, and $\Kinf=\Fq((1/T))$ is the completion of $\Fq(T)$ at the infinite place, with absolute value $|x|=q^{\deg x}$ (and $|0|=0$). Every $x\in\Kinf$ can be written uniquely as $x=[x]+\{x\}$ with $[x]\in\A$ and $|\{x\}|<1$; we put $\nrm{x}=|\{x\}|$. This is the distance from $x$ to $\A$, and it satisfies the ultrametric inequality $\nrm{x+y}\le\max(\nrm x,\nrm y)$. We write $\M_n$ for the set of monic polynomials of degree $n$, and $\PP_n\subset\M_n$ for the subset of irreducible ones. Following the usual practice, we call irreducible monic polynomials \emph{primes} and denote them by $\pi$ in the statements; inside the proofs we use the letter $P$ only when no confusion with the numerators $P_i$ of convergents can arise. For $\alpha\in\Kinf\setminus\Fq(T)$ we put
\[
\theta(\alpha)=\sup\bigl\{\theta\ge0:\ \nrm{\alpha\pi}<|\pi|^{-\theta}\text{ for infinitely many primes }\pi\bigr\}.
\]
Since $\nrm{cx}=\nrm{x}$ for $c\in\Fq^\times$, there is no loss in restricting to monic $\pi$. Since $\theta(\alpha)$ is a supremum, it does not change if the strict inequality in its definition is replaced by $\nrm{\alpha\pi}\le|\pi|^{-\theta}$; we use this freely in Sections~\ref{sec:weil} and~\ref{sec:sqfree}. The analogue of Theorem~\ref{thm:dirichlet} holds in $\Kinf$ as well (Lemma~\ref{lem:cf} below), so the question is exactly the same as over $\mathbb R$: how much of the exponent $1$ survives when the denominators are required to be irreducible?

\subsection{Diophantine sets}
\chg{Over $\A$ the problem has a linear structure} which has no exact analogue over $\mathbb Z$. For integers $n\ge1$ and $\kappa\ge0$, consider
\[
\V_n(\alpha,\kappa)=\{f\in\M_n:\ \nrm{\alpha f}<q^{-\kappa}\}.
\]
The condition $\nrm{\alpha f}<q^{-\kappa}$ says that the coefficients of $T^{-1},\dots,T^{-\kappa}$ in $\alpha f$ vanish, and each of these coefficients is an affine linear function of the coefficients of $f$. So $\V_n(\alpha,\kappa)$ is an affine $\Fq$-subspace of $\M_n\cong\Fq^n$. We call it a \emph{Diophantine set}. If it has the expected size $q^{n-\kappa}$ and if the primes are equidistributed in it, it contains about $q^{n-\kappa}/n$ primes. Thus $\theta(\alpha)\ge\theta$ as soon as $\V_n(\alpha,\lceil\theta n\rceil)$ contains a prime for infinitely many $n$. Heuristically this should be the case for every $\theta<1$; the difficulty is to prove it. It is natural to conjecture that $\theta(\alpha)\ge1$ for every $\alpha\in\Kinf\setminus\Fq(T)$, with equality when $\alpha$ has bounded partial quotients, since then $\nrm{\alpha f}\gg|f|^{-1}$ for all nonzero $f\in\A$. By Proposition~\ref{prop:level} with $\nu=0$, such a bound would follow from an asymptotic formula for primes in a single arithmetic progression whose modulus is almost as large as the primes themselves, that is, from a function field analogue of conjectures of Linnik and Montgomery type. Over $\Fq(T)$, results of this type were proved by Bank, Bary-Soroker and Rosenzweig~\cite{BBR} in the limit $q\to\infty$. For our problem, however, one needs them for fixed $q$, which seems to be much harder.

\subsection{Continued fractions}
Let $\alpha=[a_0;a_1,a_2,\dots]$ be the continued fraction expansion of $\alpha$, with $a_i\in\A$ and $\deg a_i\ge1$ for $i\ge1$, and let $p_i/q_i$ be its convergents, normalised so that $q_i$ is monic. (The letter $q$ is thus used both for the size of the field and for convergent denominators; the meaning will always be clear from the context.) Put $d_i=\deg q_i$; thus $d_0=0$ and $d_{i+1}=d_i+\deg a_{i+1}$. With the convention $d_{-1}=-\infty$, the numbers $d_i+d_{i+1}$ ($i\ge-1$) increase strictly from $-\infty$ to $\infty$. Hence for every integer $X\ge0$ there is a unique index $j=j_\alpha(X)\ge-1$ with
\begin{equation}\label{eq:jdef}
d_j+d_{j+1}\le X<d_{j+1}+d_{j+2}.
\end{equation}
We measure the size of the partial quotients by
\begin{align*}
A_\alpha(X)&=\max\{\deg a_{i+1}:\ i\ge0,\ d_i\le X\},\\
E_\alpha(X)&=\max\{\deg a_i:\ i\in\{j+1,j+2\},\ i\ge1\},\qquad j=j_\alpha(X),\\
\tau(\alpha)&=\limsup_{i\to\infty}\frac{\deg a_{i+1}}{d_i},\qquad
\sigma(\alpha)=\liminf_{j\to\infty}\frac{\max(\deg a_{j+1},\deg a_{j+2})}{d_j+d_{j+1}} .
\end{align*} 
Here $A_\alpha$ is cumulative, whereas $E_\alpha(X)\le A_\alpha(X)$ only sees the two partial quotients at the scale $X$. When $j_\alpha(X)=-1$ the index set in the definition of $E_\alpha(X)$ is $\{1\}$, so that $E_\alpha(X)=\deg a_1$; in all cases $E_\alpha(X)\ge1$. The irrationality exponent of $\alpha$ is $2+\tau(\alpha)$ (Lemma~\ref{lem:cf}). If the partial quotients of $\alpha$ have bounded degree, then $\tau(\alpha)=\sigma(\alpha)=0$; this is the case, for example, for quadratic irrationals, whose expansions are eventually periodic. The condition $\sigma(\alpha)=0$ is much weaker than $\tau(\alpha)=0$: it only asks for infinitely many pairs of consecutive partial quotients which are small compared with the convergent denominators.

\subsection{Main results}\label{sec:mainresults}
Our first two theorems give an asymptotic formula for primes in a Diophantine set when $\kappa<n/2$. In both, the partial quotients of $\alpha$ enter only through the single quantity $E_\alpha(n+\kappa)$. They differ in the condition on $q$ and in how strongly $E_\alpha(n+\kappa)$ enters (see the comparison after Theorem~\ref{thm:main2}). As usual, $\Lambda(f)=\deg\pi$ if $f$ is a power of a prime $\pi$, and $\Lambda(f)=0$ otherwise.

\begin{theorem}\label{thm:main}
Let $p$ be an odd prime, $0<\varepsilon\le\frac1{10}$, and let $q$ be a power of $p$ with
\begin{equation}\label{eq:qcond}
q>p^2e^2\Bigl(1+\frac2\varepsilon\Bigr)^2 .
\end{equation}
There are $c=c(q,\varepsilon)>0$ and $n_0=n_0(q,\varepsilon)$ such that the following holds. For every $\alpha\in\Kinf\setminus\Fq(T)$, every $n\ge n_0$ and every integer $\kappa\ge1$ with
\begin{equation}\label{eq:range}
\kappa+\tfrac12E_\alpha(n+\kappa)\le\bigl(\tfrac12-3\varepsilon\bigr)n,
\end{equation}
we have
\[
\sum_{f\in\V_n(\alpha,\kappa)}\Lambda(f)=q^{n-\kappa}\bigl(1+O(q^{-cn})\bigr),
\]
where the implied constant depends only on $q$ and $\varepsilon$.
\end{theorem}

\begin{theorem}\label{thm:main2}
Let $p$ be an odd prime, $0<\eta\le\frac1{12}$, and let $q$ be a power of $p$ with
\begin{equation}\label{eq:qcond2}
q>4p^2e^2\Bigl(1+\frac2\eta\Bigr)^2 .
\end{equation}
There are $c=c(q,\eta)>0$ and $n_0=n_0(q,\eta)$ such that the following holds. For every $\alpha\in\Kinf\setminus\Fq(T)$, every $n\ge n_0$ and every integer $\kappa$ with $\kappa\ge\frac n3$ and
\begin{equation}\label{eq:range2}
2\kappa+\tfrac12E_\alpha(n+\kappa)\le(1-4\eta)\,n,
\end{equation}
we have
\[
\sum_{f\in\V_n(\alpha,\kappa)}\Lambda(f)=q^{n-\kappa}\bigl(1+O(q^{-cn})\bigr),
\]
where the implied constant depends only on $q$ and $\eta$.
\end{theorem}

Since $\kappa\ge\frac n3$ and $E_\alpha\ge1$, condition~\eqref{eq:range2} can only hold when $\eta<\frac1{12}$, so the restriction $\eta\le\frac1{12}$ in Theorem~\ref{thm:main2} is no loss.

\emph{Comparison.} Condition~\eqref{eq:range} is equivalent to $2\kappa+E_\alpha(n+\kappa)\le(1-6\varepsilon)n$, whereas~\eqref{eq:range2} only asks for $2\kappa+\frac12E_\alpha(n+\kappa)\le(1-4\eta)n$. Thus Theorem~\ref{thm:main2} tolerates partial quotients twice as large at the scale $n+\kappa$, at the price of the requirement $\kappa\ge\frac n3$ and of the factor $4$ in~\eqref{eq:qcond2} compared with~\eqref{eq:qcond}. This is why both theorems are kept: Theorem~\ref{thm:main} gives the smaller values of $q$ in Table~\ref{tab:q} for $\theta$ close to $\frac13$, and Theorem~\ref{thm:main2} gives the larger range of $\theta$.

\emph{Uniformity.} In Theorems~\ref{thm:main} and~\ref{thm:main2} the constants $c$ and $n_0$ and the implied constants depend only on $(q,\varepsilon)$, respectively $(q,\eta)$, and not on $\alpha$; the element $\alpha$ enters only through the hypotheses~\eqref{eq:range} and~\eqref{eq:range2}. In the proofs, apart from factors polynomial in $n$, every factor $q^{o(n)}$ comes from Lemma~\ref{lem:elem}(d) (divisor bounds for polynomials of degree at most $n+\kappa$) or from Lemma~\ref{lem:elem}(f), and both are uniform in all polynomials involved; Lemma~\ref{lem:count} is uniform in $\alpha$; and the M\"obius bounds of Section~\ref{sec:SS} are uniform in the coefficients of the progressions.

Under either hypothesis $|\V_n(\alpha,\kappa)|=q^{n-\kappa}$ (see Sections~\ref{sec:proof} and~\ref{sec:proof2}), so these theorems say that $\Lambda$ has its expected mean value on $\V_n$. \chg{By Lemma~\ref{lem:elem}(e), prime powers contribute at most $5q^{n/2}$}, and $n/2<n-\kappa-cn$, so they give
\begin{equation}\label{eq:primecount}
\#\{\pi\in\PP_n:\ \nrm{\alpha\pi}<q^{-\kappa}\}=\frac{q^{n-\kappa}}{n}\bigl(1+O(q^{-cn})\bigr).
\end{equation}

\begin{corollary}\label{cor:tau}
Let $p,\varepsilon,q$ be as in Theorem~\ref{thm:main}.
\begin{enumerate}[label=\rm(\alph*)]
\item Let $\tau\ge0$ and $\theta>0$ with $\theta(2+\tau)<1-\tau-6\varepsilon$. If $\tau(\alpha)<\tau$, then~\eqref{eq:primecount} holds with $\kappa=\lceil\theta n\rceil$ for all large $n$, and $\theta(\alpha)\ge\theta$. In particular, if $\tau(\alpha)=0$ then $\theta(\alpha)\ge\frac12-3\varepsilon$.
\item Let $\sigma\ge0$ and $\theta>0$ with $\theta(2+\sigma)<1-\sigma-6\varepsilon$. If $\sigma(\alpha)<\sigma$, then~\eqref{eq:primecount} holds with $\kappa=\lceil\theta n\rceil$ for infinitely many $n$, and $\theta(\alpha)\ge\theta$.
\end{enumerate}
\end{corollary}

\begin{corollary}\label{cor:sigma}
Let $p,\eta,q$ be as in Theorem~\ref{thm:main2}, let $\sigma\ge0$, and let $\frac13\le\theta<\frac12$ with
\[
\theta(4+\sigma)<2-\sigma-8\eta .
\]
If $\sigma(\alpha)<\sigma$, then~\eqref{eq:primecount} holds with $\kappa=\lceil\theta n\rceil$ for infinitely many $n$, and $\theta(\alpha)\ge\theta$. In particular, if $\sigma(\alpha)=0$ then $\theta(\alpha)\ge\frac12-2\eta$.
\end{corollary}

For $\sigma>0$ the exponent $\frac{1-\sigma}{2+\sigma}$ of Corollary~\ref{cor:tau}(b) is smaller than the exponent $\frac{2-\sigma}{4+\sigma}$ of Corollary~\ref{cor:sigma}, since $(1-\sigma)(4+\sigma)=4-3\sigma-\sigma^2<(2-\sigma)(2+\sigma)$; but Corollary~\ref{cor:tau}(b) needs only~\eqref{eq:qcond} and no lower bound on $\theta$.

For well approximable $\alpha$ we use Weil's theorem directly. The following result holds for every prime power $q$; part~(a) is the analogue of Matom\"aki's theorem, proved in~\cite{BM} for odd $q\ge7$. The restriction on $q$ in~\cite{BM} came from the final comparison of error terms there\chg{, namely the inequality $2^{2/3}q\,(q-1)^{-2/3}<q^{2/3}$ at the end of the proof of~\cite[Theorem~3.3]{BM}, which is equivalent to $(q-1)^2>4q$ and so holds only for $q\ge7$. It arises because the number of divisors of the modulus $f$ was bounded by $2^{\deg f}$, and $1/\varphi(f)$ by $(q-1)^{-\deg f}$.} In Section~\ref{sec:weil} the same quantities are controlled by Lemmas~\ref{lem:elem}(d) and~\ref{lem:Wstar}, which give losses of size $2^{\omega(f)}=q^{o(\deg f)}$ and $|f|/\varphi(f)=q^{o(\deg f)}$ only; together with the exact reduction to residue classes in Lemma~\ref{lem:reduction}, this is why Theorem~\ref{thm:weil}(a) holds for every $q$.

\begin{theorem}\label{thm:weil}
Let $q$ be any prime power and $\alpha\in\Kinf\setminus\Fq(T)$.
\begin{enumerate}[label=\rm(\alph*)]
\item $\theta(\alpha)\ge\frac13$.
\item If $\tau(\alpha)>0$ and $\tau=\min(\tau(\alpha),1)$, then $\theta(\alpha)\ge\dfrac{1+\tau}{3+\tau}$.
\end{enumerate}
\end{theorem}

Combining the two methods gives an exponent beyond $\frac13$ for every $\alpha$.

\begin{corollary}\label{cor:all}
Let $p$ be an odd prime and $\frac13\le\theta<\sqrt2-1=0.41421\ldots$. There is $q_0(p,\theta)$ such that $\theta(\alpha)\ge\theta$ for every power $q>q_0(p,\theta)$ of $p$ and every $\alpha\in\Kinf\setminus\Fq(T)$. One can take $q_0(p,\theta)=4p^2e^2(1+2/\eta)^2$ for any $\eta\in(0,\frac1{12}]$ with
\[
8\eta<2-\sigma_1-\theta(4+\sigma_1),\qquad \sigma_1=\frac{\tau_1(1+\tau_1)}{2+\tau_1},\qquad \tau_1\in(0,1),\ \ \tau_1>\frac{3\theta-1}{1-\theta}.
\]
For $\theta<\frac{\sqrt{17}-1}8=0.39038\ldots$ one can alternatively take $q_0(p,\theta)=p^2e^2(1+2/\varepsilon)^2$ with $6\varepsilon<(1-\sigma_1)-\theta(2+\sigma_1)$, where $\sigma_1$ and $\tau_1$ are as above.
\end{corollary}

The number $\sqrt2-1$ is the value at which the exponent $\frac{1+\tau}{3+\tau}$ of Theorem~\ref{thm:weil}(b) and the exponent $\frac{2-\sigma}{4+\sigma}$ of Corollary~\ref{cor:sigma}, with $\sigma=\frac{\tau(1+\tau)}{2+\tau}$, meet; this happens at $\tau=\sqrt2-1$. Table~\ref{tab:q} lists admissible values of $q$.

\begin{table}[ht]\centering
\begin{tabular}{lccc}\toprule
$\theta$ & via Theorem~\ref{thm:main} & via Theorem~\ref{thm:main2} & smallest $q$ for $p=3,5,7$\\\midrule
$\frac13+\frac1{126}$ & $13026\,p^2$ & $21897\,p^2$ & $3^{11},\ 5^{8},\ 7^{7}$\\
$0.35$ & $18290\,p^2$ & $27128\,p^2$ & $3^{11},\ 5^{9},\ 7^{8}$\\
$0.36$ & $30418\,p^2$ & $36296\,p^2$ & $3^{12},\ 5^{9},\ 7^{8}$\\
$0.37$ & $63570\,p^2$ & $51998\,p^2$ & $3^{12},\ 5^{9},\ 7^{8}$\\
$0.38$ & $230167\,p^2$ & $82667\,p^2$ & $3^{13},\ 5^{10},\ 7^{8}$\\
$0.39$ & $1.56\cdot10^8\,p^2$ & $156980\,p^2$ & $3^{13},\ 5^{10},\ 7^{9}$\\
$0.40$ & -- & $432923\,p^2$ & $3^{14},\ 5^{11},\ 7^{9}$\\
$0.41$ & -- & $4.68\cdot10^6\,p^2$ & $3^{16},\ 5^{12},\ 7^{10}$\\\bottomrule
\end{tabular}
\caption{Sufficient lower bounds for $q$ in Corollary~\ref{cor:all} (all $\alpha$), with $\tau_1$ taken just above $\frac{3\theta-1}{1-\theta}$. The thresholds are rounded \chg{up} and the inequalities on $\varepsilon,\eta$ are strict.}\label{tab:q}
\end{table}

For $\alpha$ with $\tau(\alpha)=0$, Corollary~\ref{cor:tau}(a) gives every $\theta\in[\frac13,\frac12)$ as soon as $q>e^2\bigl(1+\frac{12}{1-2\theta}\bigr)^2p^2$ (for such $\theta$ the admissible $\varepsilon$ satisfy $\varepsilon<\frac{1-2\theta}6\le\frac1{18}$, so the restriction $\varepsilon\le\frac1{10}$ is automatic); for example, $\theta=\frac13+\frac1{126}$ needs $q>11124\,p^2$, so $q=3^{11}$, $5^8$ or $7^7$ suffices.

\subsection{Squarefree moduli}
In the classical approach the moduli are the convergent denominators $f=q_m$ of $\alpha$, and these are imposed on us by $\alpha$. When they happen to be squarefree, a different and much stronger input becomes available. Sawin~\cite{Sawin} proved an asymptotic formula for primes in a single residue class modulo a squarefree $f$, with a level of distribution which tends to $1$ as $q\to\infty$, and with a condition on $q$ that does not depend on the characteristic. For $q>e^2$ let $\lambda(q)\in(0,1)$ be defined by
\begin{equation}\label{eq:lambdaq}
q=e^2\,4^{\frac1{1-\lambda(q)}}\Bigl(\frac1{1-\lambda(q)}-\frac12\Bigr)^2 .
\end{equation}
Then $\lambda$ is increasing, $\lambda(2956)>\frac23$, $\lambda(23173)>\frac34$, and $\lambda(q)\to1$ as $q\to\infty$. (The right side of~\eqref{eq:lambdaq} equals $400e^2=2955.62\ldots$ at $\lambda=\frac23$ and $3136e^2=23172.07\ldots$ at $\lambda=\frac34$.)

\begin{theorem}\label{thm:sqfree}
Let $q>e^2$ be any prime power, and let $\alpha\in\Kinf\setminus\Fq(T)$. Suppose that, for infinitely many $m$, at least one of the $q$ polynomials $q_m+c\,q_{m-1}$ \textup{(}$c\in\Fq$\textup{)} is squarefree. Then
\[
\theta(\alpha)\ge2\lambda(q)-1 .
\]
The hypothesis holds, already with $c=0$, for all $\alpha$ with $|\alpha|<1$ outside a set of Haar measure at most $2/q$ \textup{(}with the Haar measure normalised to total mass $1$ on $\{|x|<1\}$\textup{)}, and likewise in every translate of this set by $\A$.
\end{theorem}

We now compare Theorem~\ref{thm:sqfree} with the main theorems, since the two use the moduli in different ways.
\begin{itemize}[leftmargin=2em] 
\item \emph{Arbitrary moduli.} Theorems~\ref{thm:main} and~\ref{thm:main2} need no arithmetic condition at all on the moduli. They apply for odd $p$ and $q$ large in terms of $p$, and give Corollaries~\ref{cor:tau}--\ref{cor:all}. No splitting into a squarefree and a non-squarefree case occurs in their proofs.
\item \emph{Squarefree moduli.} When $\alpha$ has infinitely many squarefree convergent denominators (or squarefree members of the pencils above), Theorem~\ref{thm:sqfree} applies to \emph{every} prime power $q$, including prime $q$. The exponent $2\lambda(q)-1$ exceeds $\frac13$ for $q>2956$, exceeds $\frac12$ for $q>23172$, and tends to $1$.
\end{itemize} 
In particular, for prime $q>2956$ we obtain an exponent beyond $\frac13$ for every $\alpha$ satisfying the hypothesis of Theorem~\ref{thm:sqfree}, and hence for all $\alpha$ outside a set of measure at most $2/q$. For general $\alpha$ and prime $q$ the problem remains open.

\subsection{The idea of the proof}
Write $S=\sum_{f\in\V_n}\Lambda(f)$ and expand $\Lambda(f)=-\sum_{D\mid f}\mu(D)\deg D$, where $D$ runs over monic divisors. This gives
\[
S=-\sum_{k=1}^{n}k\sum_{D\in\M_k}\mu(D)\,\#\{B\in\M_{n-k}:\ DB\in\V_n\}.
\]
We treat small and large divisors $D$ in different ways.

\emph{Small divisors.} For fixed $D$, the set of cofactors $B$ is again a Diophantine set, namely $\V_{n-k}(D\alpha,\kappa)$. Being an affine subspace, its size is exactly $q^{n-k-\kappa}$, unless a certain Hankel matrix drops rank (Lemma~\ref{lem:hankel}). The rank can drop only if $GD\alpha$ is unusually close to $\A$ for some nonzero $G$ of degree $<\kappa$, and the number of such $D$ is controlled by the partial quotients of $\alpha$ at the scale $n+\kappa$ (Lemma~\ref{lem:count}). So, for small $D$, there is no error term at all apart from these exceptional $D$. The main term comes from $k=1$ alone, because $\sum_{D\in\M_k}\mu(D)=0$ for $k\ge2$.

\emph{Large divisors.} Here we interchange the summations. For fixed $B$, the set $\F_B$ of admissible $D$ is the Diophantine set $\V_k(B\alpha,\kappa)$. Let $Q_j,Q_{j+1}$ be the consecutive convergent denominators of $B\alpha$ with $\deg Q_j+\deg Q_{j+1}\le k+\kappa<\deg Q_{j+1}+\deg Q_{j+2}$. Then (Lemma~\ref{lem:fibre})
\[
\F_B=\bigl\{UQ_j+WQ_{j+1}:\ |U|<q^{s-\kappa},\ W\in\M_{k-s}\bigr\},\qquad s=\deg Q_{j+1}.
\]
This description is the key to the argument. For fixed $U$, the set is an arithmetic progression modulo $Q_{j+1}$ of length $q^{k-s}$; for fixed $W$, it is a progression modulo $Q_j$ along the interval $\{|U|<q^{s-\kappa}\}$. Since $(k-s)+(s-\kappa)=k-\kappa$, at least one of the two directions is long as soon as $k$ exceeds $\kappa$ by a fixed proportion. Sawin and Shusterman proved cancellation of $\mu$ along such linear sequences, for arbitrary moduli and $q$ large in terms of $p$: along complete sets $\M_t$ in~\cite[Theorem~4.5]{SS}, and along arbitrary intervals, uniformly in the coefficients, in~\cite[Theorem~5.1]{SS2}. In Theorem~\ref{thm:main} only the first direction is used, and the remaining $B$ are counted. In Theorem~\ref{thm:main2} both directions are used, and no counting is needed for large divisors.

The two ranges overlap exactly when $\kappa$ is below $n/2$, up to losses caused by large partial quotients; this is where the restriction $\theta<\frac12$ comes from (see \S\ref{sec:remarks}). For Theorem~\ref{thm:weil} we use the classical method: the Diophantine set is a union of residue classes modulo a convergent denominator of $\alpha$, and Weil's theorem, together with an averaging over the residue classes, gives an asymptotic formula for moduli up to $X^{1-\theta}$, where $X=q^n$. The two approaches \chg{are complementary}: the first works best when the partial quotients of $\alpha$ are small, the second when they are large, and Corollary~\ref{cor:all} is obtained by balancing them.

\subsection{Relation to the level of distribution}
For a convergent $a/f$ of $\alpha$ with $\deg f\approx\frac{n+\kappa}2$, the set $\V_n$ is a union of about $q^{\deg f-\kappa}$ residue classes modulo $f$ (Lemma~\ref{lem:reduction}). The Riemann Hypothesis handles $|f|\le X^{2/3-\varepsilon}$, which corresponds to $\theta=\frac13$. Theorem~\ref{thm:main} handles moduli up to $|f|\approx X^{3/4-O(\varepsilon)}$ when $\tau(\alpha)=0$. Sawin and Shusterman~\cite[Theorem~1.9]{SS} (numbered as in the arXiv version) proved an asymptotic formula for primes in a \emph{single} progression with $|f|\le X^{1/2+\delta}$, $\delta<\frac1{126}$, for squarefree $f$, and they remark that their convolution method cannot go beyond $X^{2/3}$ for a single progression. Sawin~\cite{Sawin} obtained levels tending to $1$ as $q\to\infty$, again for squarefree moduli. A level $\lambda$ used for one modulus at a time only gives the exponent $2\lambda-1$ (Proposition~\ref{prop:level}), and it requires the modulus, a convergent denominator dictated by $\alpha$, to be squarefree. Our main argument uses instead M\"obius estimates along linear sequences, which hold for all moduli, and the averaging over the residue classes is essential for it.

\subsection{Organisation of the paper}
In Section~\ref{sec:prelim} we collect elementary facts and prove what we need on continued fractions in $\Kinf$. Section~\ref{sec:dioph} contains three lemmas on Diophantine sets, which are the main new elementary tools. Section~\ref{sec:weil} proves Theorem~\ref{thm:weil}, and Section~\ref{sec:sqfree} proves Theorem~\ref{thm:sqfree}. Section~\ref{sec:SS} states the input from~\cite{SS,SS2} and estimates M\"obius sums over fibres. Sections~\ref{sec:proof} and~\ref{sec:proof2} prove Theorems~\ref{thm:main} and~\ref{thm:main2}, and Section~\ref{sec:cor} proves the corollaries. In Section~\ref{sec:remarks} we discuss the limits of the method and explain why the exponent $\frac12$ is a natural barrier for it.

\subsection{Use of AI tools}
The author used Claude (Anthropic, Claude Opus~5.5) as a research and writing assistant. The research problem, the overall strategy and the decisions on scope are the author's. Based on the author's material and instructions, Claude observed that the proof of Theorem~\ref{thm:main} needs the partial quotients of $\alpha$ only through $E_\alpha(n+\kappa)$, derived Corollary~\ref{cor:tau}(b) and the range $\theta<\frac{\sqrt{17}-1}8$ in Corollary~\ref{cor:all}, recomputed Table~\ref{tab:q}, drafted Remarks~\ref{rem:lacunary} and~\ref{rem:sqfree}, checked the proofs, and assisted with the exposition. The author has checked all statements and proofs and takes full responsibility for the content.
 
\section{Preliminaries}\label{sec:prelim}

In this section we fix notation and collect the elementary facts used later. Everything here is standard, but we give proofs, or precise references, so that the paper can be read without consulting other sources for these points.

\subsection{Notation}
All divisors are monic. We write $\mu$, $\Lambda$, $\varphi$, $\tau_d$ and $\omega$ for the M\"obius function, the von Mangoldt function ($\Lambda(f)=\deg P$ if $f=P^r$ with $P$ prime and $r\ge1$, and $\Lambda(f)=0$ otherwise), the Euler function, the divisor function and the number of distinct prime factors. As in~\cite{SS2}, $\mu$ is defined on all nonzero polynomials: $\mu(f)=0$ if $f$ is divisible by the square of a nonconstant polynomial, and $\mu(f)=(-1)^r$ otherwise, where $r$ is the number of monic prime factors of $f$. For a real number $a$, the condition $|U|<q^{a}$ on $U\in\A$ means $\deg U<a$; it allows only $U=0$ when $a\le0$, and allows exactly $q^{a}$ polynomials when $a$ is a positive integer. The notation $o(n)$ stands for a quantity bounded by $\eta'n$ for every $\eta'>0$, once $n$ is large in terms of $\eta'$ and $q$.

\begin{lemma}\label{lem:elem}
With the notation above, we have the following. 
\begin{enumerate}[label=\rm(\alph*)]
\item $\#\M_n=q^n$, and $\#\PP_n=\frac{q^n}n+O\bigl(\frac{q^{n/2}}n\bigr)$.
\item $\sum_{D\in\M_k}\mu(D)$ equals $1$, $-q$, $0$ for $k=0$, $k=1$, $k\ge2$, respectively.
\item For $f\in\M_n$ with $n\ge1$: $\Lambda(f)=-\sum_{D\mid f}\mu(D)\deg D$.
\item For nonzero $g$ with $\deg g\ge1$: $\omega(g)=o(\deg g)$, $\tau_d(g)=q^{o(\deg g)}$, and $|g|/\varphi(g)=q^{o(\deg g)}$.
\item $\sum_{P^e\in\M_n,\ e\ge2}\deg P\le q^{n/2}+n\,q^{n/3}\le5\,q^{n/2}$, where $P$ runs over primes.
\item For $Q\ne0$ and integers $0\le t\le n$ with $\deg Q\le n$: $\bigl|\sum_{W\in\M_t}\mu(WQ)\bigr|\le q^{1+o(n)}$, where the $o(n)$ depends only on $q$, not on $Q$ or $t$.
\end{enumerate}
\end{lemma}

\begin{proof}
(a) is the prime polynomial theorem~\cite[Theorem~2.2]{Rosen}.

(b) The generating function is $\sum_{D}\mu(D)u^{\deg D}=\prod_P(1-u^{\deg P})=1/Z(u)=1-qu$, where $Z(u)=\sum_Du^{\deg D}=(1-qu)^{-1}$.

(c) Since $\sum_{D\mid f}\Lambda(D)=\deg f$, M\"obius inversion gives $\Lambda(f)=\sum_{D\mid f}\mu(D)\deg(f/D)=\deg f\sum_{D\mid f}\mu(D)-\sum_{D\mid f}\mu(D)\deg D$, and the first sum vanishes for $\deg f\ge1$.

(d) Let $N=\deg g$ and $m=\lfloor\frac12\log_qN\rfloor$. At most $\sum_{i\le m}q^i\le2q^m\le2\sqrt N$ primes have degree $\le m$, and at most $N/(m+1)$ prime factors of $g$ have degree $>m$. Hence $\omega(g)\le2\sqrt N+N/(m+1)=o(N)$. Write $g=\prod P_i^{e_i}$. The primes of degree $\le m$ contribute at most $(N+1)^{2\sqrt N}$ to $\tau_d(g)=\prod(e_i+1)$. For the others, $\sum e_i\le N/(m+1)$ and $e+1\le2^e$, so they contribute at most $2^{N/(m+1)}$. Both are $q^{o(N)}$. Finally $|g|/\varphi(g)=\prod_{P\mid g}(1-|P|^{-1})^{-1}\le(1-q^{-1})^{-\omega(g)}=q^{o(N)}$.

(e) Since $\sum_{d\mid m}d\,\#\PP_d=q^m$, we have $\#\PP_d\le q^d/d$. For each $e\ge2$ with $e\mid n$, the prime powers $P^e\in\M_n$ have $\deg P=n/e$, so they contribute at most $\frac ne\#\PP_{n/e}\le q^{n/e}$. Summing over $2\le e\le n$ gives at most $q^{n/2}+nq^{n/3}$, and $nq^{n/3}\le n2^{-n/6}q^{n/2}\le4\,q^{n/2}$, since $\max_{n\ge1}n2^{-n/6}<3.2$.

(f) We have $\mu(WQ)=\mu(W)\mu(Q)$ if $(W,Q)=1$, and $\mu(WQ)=0$ otherwise. Moreover
\[
\sum_{(W,Q)=1}\mu(W)u^{\deg W}=\prod_{P\nmid Q}(1-u^{\deg P})=(1-qu)\prod_{P\mid Q}(1-u^{\deg P})^{-1}.
\]
The coefficient of $u^t$ on the right is $c_t-qc_{t-1}$ (with $c_{-1}=0$), where $c_m\ge0$ is the number of $W\in\M_m$ all of whose prime factors divide $Q$. Since $\mu(WQ)=0$ for all $W$ when $Q$ is not squarefree, and $|\mu(Q)|\le1$, we get
\begin{equation}\label{eq:cbound}
\Bigl|\sum_{W\in\M_t}\mu(WQ)\Bigr|\le|c_t-qc_{t-1}|\le q\max(c_t,c_{t-1}).
\end{equation}
It remains to show that $c_m\le q^{o(n)}$ for $0\le m\le n$, uniformly in $Q$ with $\deg Q\le n$.

First we bound $\omega(Q)$ uniformly. Put $m_0=\lfloor\frac12\log_qn\rfloor$. As in the proof of (d), at most $2q^{m_0}\le2\sqrt n$ primes have degree $\le m_0$, and at most $\deg Q/(m_0+1)\le n/(m_0+1)$ prime factors of $Q$ have degree $>m_0$. Hence
\[
\omega(Q)\le w(n):=2\sqrt n+\frac{n}{m_0+1},\qquad w(n)=o(n),
\]
and $w(n)$ depends only on $n$ and $q$.

Next we use Rankin's trick. Fix $0<\sigma<1$. Every $W$ counted by $c_m$ satisfies $(q^m/|W|)^\sigma=1$, and all other terms below are positive, so
\[
c_m\le\sum_{\substack{W\ \mathrm{monic}\\ P\mid W\Rightarrow P\mid Q}}\Bigl(\frac{q^{m}}{|W|}\Bigr)^{\sigma}
=q^{\sigma m}\prod_{P\mid Q}\bigl(1-|P|^{-\sigma}\bigr)^{-1}
\le q^{\sigma n}\bigl(1-q^{-\sigma}\bigr)^{-w(n)} .
\]
Now let $\varepsilon'>0$ and take $\sigma=\varepsilon'/2$. Since the base $(1-q^{-\sigma})^{-1}$ is a constant depending only on $q$ and $\varepsilon'$, and $w(n)=o(n)$, we have $(1-q^{-\sigma})^{-w(n)}\le q^{\varepsilon'n/2}$ for all $n\ge n_0(q,\varepsilon')$. Hence $c_m\le q^{\varepsilon'n}$ for $n\ge n_0(q,\varepsilon')$, all $0\le m\le n$ and all $Q$ with $\deg Q\le n$; that is, $c_m\le q^{o(n)}$ uniformly. Together with~\eqref{eq:cbound} this proves (f).

We remark that the cruder bound $c_m\le(m+1)^{\omega(Q)}$, obtained by counting exponent vectors, is not sufficient here. For example, if $Q$ is the product of all primes of degree at most $m_1$, then $\deg Q\asymp q^{m_1}$ and $\omega(Q)\asymp\deg Q/m_1$, so that $(\deg Q+1)^{\omega(Q)}$ is of size roughly $q^{\deg Q}$ (for instance about $q^{1.04\deg Q}$ when $q=7$, $m_1=7$). The point of Rankin's trick is that it replaces the growing base $m+1$ by the constant base $(1-q^{-\sigma})^{-1}$.
\end{proof}

\subsection{Dirichlet characters}

\begin{lemma}[Weil]\label{lem:PNT}
Let $f\in\A$ be monic of degree $F\ge1$, $\chi$ a Dirichlet character modulo $f$, and $N\ge1$. Then
\[
\sum_{\substack{P\in\PP_N\\ P\nmid f}}\chi(P)=\mathbf 1_{\chi=\chi_0}\frac{q^N}{N}+O\Bigl(\frac{(F+1)\,q^{N/2}}{N}\Bigr),
\]
with an absolute implied constant.
\end{lemma}

\begin{proof}
Let $\chi\ne\chi_0$. The $L$-function $L(u,\chi)=\sum_{g}\chi(g)u^{\deg g}$, the sum over monic $g$, is a polynomial in $u$ of degree at most $F-1$, because $\sum_{g\in\M_m}\chi(g)=0$ for $m\ge F$. By Weil's theorem~\cite[Ch.~4 and Ch.~9]{Rosen}, its inverse roots $\gamma_1,\dots,\gamma_r$ ($r\le F-1$) have absolute value $\sqrt q$ or $1$. Taking the logarithmic derivative of the Euler product gives
\[
\sum_{g\in\M_N}\Lambda(g)\chi(g)=-\sum_{i=1}^r\gamma_i^N,
\]
which has absolute value at most $(F-1)q^{N/2}$. The prime powers $P^e$ with $e\ge2$ contribute at most $5q^{N/2}$ to the left side by Lemma~\ref{lem:elem}(e), and the primes contribute $N\sum_{P\in\PP_N}\chi(P)$. Hence $|\sum_{P\in\PP_N}\chi(P)|\le(F+4)q^{N/2}/N$. For $\chi=\chi_0$, use Lemma~\ref{lem:elem}(a), and note that at most $F/N$ primes of degree $N$ divide $f$.
\end{proof}

\subsection{Continued fractions}
We recall the facts we need, with proofs for completeness; see~\cite{Lasjaunias,Schmidt} for more. Let $\beta\in\Kinf\setminus\Fq(T)$. Its continued fraction expansion $\beta=[b_0;b_1,b_2,\dots]$ is obtained by $\beta_0=\beta$, $b_i=[\beta_i]$, $\beta_{i+1}=1/(\beta_i-b_i)$. Thus $b_i\in\A$, $\deg b_i\ge1$ for $i\ge1$, and $|\beta_{i}|=|b_{i}|$ for $i\ge1$. Define $P_{-1}=1$, $Q_{-1}=0$, $P_0=b_0$, $Q_0=1$ and
\[
P_{i+1}=b_{i+1}P_i+P_{i-1},\qquad Q_{i+1}=b_{i+1}Q_i+Q_{i-1}\qquad(i\ge0).
\]

\begin{lemma}\label{lem:cf}
Put $d_i=\deg Q_i$ and $\delta_i=Q_i\beta-P_i$, with $d_{-1}=-\infty$. Then for all $i\ge-1$:
\begin{enumerate}[label=\rm(\alph*)]
\item $d_{i+1}=d_i+\deg b_{i+1}$ for $i\ge0$, and $Q_{i+1}P_i-Q_iP_{i+1}=(-1)^{i+1}$;
\item $\beta=\dfrac{\beta_{i+1}P_i+P_{i-1}}{\beta_{i+1}Q_i+Q_{i-1}}$ for $i\ge0$;
\item $|\delta_i|=q^{-d_{i+1}}$, and $|\beta-P_i/Q_i|=q^{-d_i-d_{i+1}}$ for $i\ge0$.
\end{enumerate}
In particular the irrationality exponent of $\beta$ is $2+\tau(\beta)$, where $\tau(\beta)=\limsup\deg b_{i+1}/d_i$.
\end{lemma}

\begin{proof}
(a) The degree statement follows by induction from the recursion, since $\deg(b_{i+1}Q_i)>\deg Q_{i-1}$. The determinant satisfies $Q_{i+1}P_i-Q_iP_{i+1}=-(Q_iP_{i-1}-Q_{i-1}P_i)$, and equals $Q_0P_{-1}-Q_{-1}P_0=1$ for $i=-1$.

(b) By induction on $i$: for $i=0$ this is $\beta=b_0+1/\beta_1$, and the step uses $\beta_{i+1}=b_{i+1}+1/\beta_{i+2}$ and the recursion.

(c) For $i=-1$, $|\delta_{-1}|=1=q^{-d_0}$. For $i\ge0$, by (b) and (a),
\[
\delta_i=Q_i\beta-P_i=\frac{Q_iP_{i-1}-P_iQ_{i-1}}{\beta_{i+1}Q_i+Q_{i-1}}=\frac{\pm1}{\beta_{i+1}Q_i+Q_{i-1}} .
\]
Since $|\beta_{i+1}Q_i|=q^{\deg b_{i+1}+d_i}=q^{d_{i+1}}>|Q_{i-1}|$, we get $|\delta_i|=q^{-d_{i+1}}$. The second formula follows by dividing by $|Q_i|$. For the last claim, see Lemma~\ref{lem:legendre} below.
\end{proof}

From now on we normalise the convergents by dividing $P_i$ and $Q_i$ by the leading coefficient of $Q_i$, so that $Q_i$ is monic. This changes neither $|\delta_i|$ nor the fractions $P_i/Q_i$, and Lemma~\ref{lem:cf}(a) becomes
\begin{equation}\label{eq:cf}
c_i:=Q_{i+1}P_i-Q_iP_{i+1}\in\Fq^\times,\qquad|\delta_i|=q^{-d_{i+1}}\qquad(i\ge-1),
\end{equation}
with $c_{-1}=1$ and $\delta_{-1}=-1$. The index $j_\beta(X)$ is defined by~\eqref{eq:jdef}.

\begin{lemma}\label{lem:coord}
Let $i\ge-1$. Every pair $(g,R)\in\A^2$ can be written uniquely as $(g,R)=U(Q_i,P_i)+W(Q_{i+1},P_{i+1})$ with $U,W\in\A$. If $\xi=g\beta-R$, then
\[
U=c_i^{-1}\bigl(g\delta_{i+1}-Q_{i+1}\xi\bigr),\qquad W=c_i^{-1}\bigl(Q_i\xi-g\delta_i\bigr).
\]
\end{lemma}

\begin{proof}
The matrix $\begin{psmallmatrix}Q_i&Q_{i+1}\\P_i&P_{i+1}\end{psmallmatrix}$ has determinant $-c_i\in\Fq^\times$ by \eqref{eq:cf}, so it lies in $\mathrm{GL}_2(\A)$. This gives existence and uniqueness of $U,W\in\A$. Next, $\xi=g\beta-R=U(Q_i\beta-P_i)+W(Q_{i+1}\beta-P_{i+1})=U\delta_i+W\delta_{i+1}$. The system
\[
g=UQ_i+WQ_{i+1},\qquad \xi=U\delta_i+W\delta_{i+1}
\]
has determinant $Q_i\delta_{i+1}-Q_{i+1}\delta_i=Q_{i+1}P_i-Q_iP_{i+1}=c_i$, and Cramer's rule gives the formulas.
\end{proof}

\begin{lemma}\label{lem:legendre}
Let $r,Q\in\A$ be coprime, $Q\ne0$ monic, with $|\beta-r/Q|<|Q|^{-2}$. Then $r/Q$ is a convergent of $\beta$.
\end{lemma}

\begin{proof}
Choose $i\ge0$ with $d_i\le\deg Q<d_{i+1}$, and write $(Q,r)=U(Q_i,P_i)+W(Q_{i+1},P_{i+1})$ as in Lemma~\ref{lem:coord}. Here $\xi=Q\beta-r$ has $|\xi|<|Q|^{-1}$. Then $|Q_i\xi|<q^{d_i-\deg Q}\le1$ and $|Q\delta_i|=q^{\deg Q-d_{i+1}}<1$, so $|W|<1$, i.e.\ $W=0$. Hence $(Q,r)=U(Q_i,P_i)$. Since $Q$ and $r$ are coprime, $U$ is a nonzero constant, and $r/Q=P_i/Q_i$.
\end{proof}

By Lemma~\ref{lem:cf}(c), the convergent $P_i/Q_i$ satisfies $|\beta-P_i/Q_i|=|Q_i|^{-2-\deg b_{i+1}/d_i}$. Together with Lemma~\ref{lem:legendre}, this shows that the irrationality exponent of $\beta$ is $2+\tau(\beta)$.

\section{Diophantine sets}\label{sec:dioph}

In this section $\beta\in\Kinf$ and $\alpha\in\Kinf\setminus\Fq(T)$ are arbitrary. We prove three lemmas: an exact formula for the size of a Diophantine set, a bound for the number of polynomials $g$ with $g\alpha$ unusually close to $\A$, and the structure of a Diophantine set as a union of progressions.

\begin{lemma}\label{lem:hankel}
Let $\beta\in\Kinf$ and $N\ge\kappa\ge0$. Put
\[
\G=\{G\in\A:\ \deg G<\kappa,\ \nrm{G\beta}<q^{-N}\},
\]
an $\Fq$-vector space of dimension $r$. Then $\#\{B\in\M_N:\ \nrm{B\beta}<q^{-\kappa}\}$ is either $0$ or $q^{N-\kappa+r}$. In particular it equals $q^{N-\kappa}$ if $\G=\{0\}$, and in all cases
\[
\bigl|\#\{B\in\M_N:\ \nrm{B\beta}<q^{-\kappa}\}-q^{N-\kappa}\bigr|\le q^{N-\kappa}\,\#(\G\setminus\{0\}).
\]
\end{lemma}

\begin{proof}
For $x\in\Kinf$, the condition $\nrm{x}<q^{-\kappa}$ means that the coefficients of $T^{-1},\dots,T^{-\kappa}$ in $x$ vanish. Write $\{\beta\}=\sum_{\chg{m}\ge1}\chg{h_m}T^{-\chg{m}}$. For $l\ge0$, the polynomial $T^l[\beta]$ contributes no negative powers of $T$, and the coefficient of $T^{-i}$ in $T^l\{\beta\}$ is $\chg{h_{i+l}}$. Write $B=T^N+\sum_{0\le l<N}b_lT^l$ with $b=(b_l)\in\Fq^N$. Then $\nrm{B\beta}<q^{-\kappa}$ is the affine system
\[
Hb=-v,\qquad H=(\chg{h_{i+l}})_{1\le i\le\kappa,\ 0\le l<N},\qquad v=(\chg{h_{i+N}})_{1\le i\le\kappa}.
\]
Its solution set is empty or an affine space with $q^{N-\operatorname{rank}H}$ elements.

To compute the rank, identify $y\in\Fq^\kappa$ with $G=\sum_{i=1}^\kappa y_iT^{i-1}$; this is a bijection onto $\{G:\deg G<\kappa\}$. The coefficient of $T^{-(l+1)}$ in $G\{\beta\}$ is $\sum_iy_i\chg{h_{i+l}}$, the $l$-th entry of $y^{\mathsf T}H$. So $y^{\mathsf T}H=0$ iff the coefficients of $T^{-1},\dots,T^{-N}$ in $G\beta$ vanish, i.e.\ iff $\nrm{G\beta}<q^{-N}$. Hence the left kernel of $H$ is isomorphic to $\G$, and $\operatorname{rank}H=\kappa-r$.

If $r=0$, $H$ has full row rank $\kappa\le N$, so the system is consistent and has exactly $q^{N-\kappa}$ solutions. If $r\ge1$, the count is $0$ or $q^{N-\kappa+r}$, and in both cases it differs from $q^{N-\kappa}$ by at most $q^{N-\kappa}(q^r-1)=q^{N-\kappa}\#(\G\setminus\{0\})$.
\end{proof}

\begin{lemma}\label{lem:count}
Let $\alpha\in\Kinf\setminus\Fq(T)$ and let $L,M\ge0$ be integers. Put $E=E_\alpha(L+M)$ and
\[
\mathcal N(L,M)=\#\{g\in\A\setminus\{0\}:\ \deg g<L,\ \nrm{g\alpha}<q^{-M}\}.
\]
\begin{enumerate}[label=\rm(\roman*)]
\item If $L\ge M$, then $\mathcal N(L,M)\le q^{L-M+E/2}$.
\item If $L<M$, then $\mathcal N(L,M)=0$ when $M-L\ge E$, and $\mathcal N(L,M)\le q^{(E-(M-L))/2}$ otherwise.
\end{enumerate}
Since $E_\alpha\le A_\alpha$, the same holds with $E=A_\alpha(L+M)$.
\end{lemma}

\begin{proof}
Use the convergents of $\alpha$, let $j=j_\alpha(L+M)$ and $s=d_{j+1}$. Let $g$ be counted (or $g=0$), let $R$ be the polynomial nearest to $g\alpha$, and $\xi=g\alpha-R$, so that $|\xi|=\nrm{g\alpha}<q^{-M}$. Write $(g,R)=U(Q_j,P_j)+W(Q_{j+1},P_{j+1})$ as in Lemma~\ref{lem:coord}. Since $|g|<q^L$:
\begin{itemize}[leftmargin=2em]
\item $|g\delta_{j+1}|<q^{L-d_{j+2}}\le q^{s-M}$, because $L+M<s+d_{j+2}$ by~\eqref{eq:jdef}; and $|Q_{j+1}\xi|<q^{s-M}$. Hence $|U|<q^{s-M}$.
\item $|g\delta_j|<q^{L-s}$; and $|Q_j\xi|<q^{d_j-M}\le q^{L-s}$, because $d_j+s\le L+M$. Hence $|W|<q^{L-s}$.
\end{itemize}
The map $g\mapsto(U,W)$ is injective. So the number of such $g$, including $g=0$, is at most $q^{\rho}$ with $\rho=(s-M)^++(L-s)^+$.

Put $X=L+M$. Since $d_j+s\le X<s+d_{j+2}$, we have
\begin{equation}\label{eq:halfE}
s-\tfrac X2\le\tfrac{s-d_j}2=\tfrac12\deg a_{j+1}\le\tfrac E2,\qquad \tfrac X2-s<\tfrac{d_{j+2}-s}2=\tfrac12\deg a_{j+2}\le\tfrac E2 .
\end{equation}
Both partial quotients enter the definition of $E=E_\alpha(X)$. If $j=-1$, then $s=0$; the first inequality is not available, but it is not needed below, since then $s>L$ and $s>M$ are impossible.

\emph{Case $L\ge M$.} If $M\le s\le L$, then $\rho=L-M$. If $s>L$, then $\rho=s-M=(s-\frac X2)+\frac{L-M}2\le\frac E2+\frac{L-M}2$. If $s<M$, then $\rho=L-s=(\frac X2-s)+\frac{L-M}2<\frac E2+\frac{L-M}2$. In all cases $\rho\le L-M+\frac E2$, which proves~(i).

\emph{Case $L<M$.} If $L\le s\le M$, then $\rho=0$, so only $g=0$ occurs. If $s>M$, then $\rho=s-M=(s-\frac X2)-\frac{M-L}2\le\frac{E-(M-L)}2$. If $s<L$, then $\rho=L-s=(\frac X2-s)-\frac{M-L}2<\frac{E-(M-L)}2$. In the last two subcases $\rho>0$, which is only possible if $M-L<E$, and then $\mathcal N\le q^{\rho}\le q^{(E-(M-L))/2}$. This proves~(ii).
\end{proof}

\begin{lemma}[Fibres]\label{lem:fibre}
Let $\beta\in\Kinf\setminus\Fq(T)$, $k\ge1$ and $\kappa\ge0$. Let $j=j_\beta(k+\kappa)$ and $s=d_{j+1}$.
\begin{enumerate}[label=\rm(\alph*)]
\item If $j=-1$, then $\nrm{D\beta}<q^{-\kappa}$ for all $D\in\M_k$.
\item If $j\ge0$, then
\[
\{D\in\M_k:\ \nrm{D\beta}<q^{-\kappa}\}=\bigl\{UQ_j+WQ_{j+1}:\ |U|<q^{s-\kappa},\ W\in\M_{k-s}\bigr\},
\]
and the pair $(U,W)$ is unique. The set is empty if $s>k$. If $s\le\kappa$, only $U=0$ occurs, and the set is $Q_{j+1}\M_{k-s}$. If $\kappa\le s\le k$, it has exactly $q^{k-\kappa}$ elements.
\end{enumerate}
\end{lemma}

\begin{proof}
(a) If $j=-1$, then $k+\kappa<d_0+d_1=d_1$. For $D\in\M_k$, $D\beta-DP_0=D\delta_0$, so $\nrm{D\beta}\le|D\delta_0|=q^{k-d_1}<q^{-\kappa}$.

(b) Let $D$ be in the left side, let $R$ be the polynomial nearest to $D\beta$, and $\xi=D\beta-R$, so $|\xi|<q^{-\kappa}$. Write $(D,R)$ as in Lemma~\ref{lem:coord} with $i=j$.
\begin{itemize}[leftmargin=2em]
\item $|D\delta_{j+1}|=q^{k-d_{j+2}}<q^{s-\kappa}$, since $k+\kappa<s+d_{j+2}$; and $|Q_{j+1}\xi|<q^{s-\kappa}$. Hence $|U|<q^{s-\kappa}$.
\item $|D\delta_j|=q^{k-s}$ exactly, and $|Q_j\xi|<q^{d_j-\kappa}\le q^{k-s}$, since $d_j+s\le k+\kappa$. By the ultrametric inequality $|W|=q^{k-s}$.
\end{itemize}
In particular $s\le k$, since $W\ne0$ is a polynomial. Moreover $\deg(UQ_j)<s-\kappa+d_j\le k$, so the leading term of $D$ equals that of $WQ_{j+1}$. As $D$ and $Q_{j+1}$ are monic, $W$ is monic of degree $k-s$.

Conversely, let $D=UQ_j+WQ_{j+1}$ with $|U|<q^{s-\kappa}$ and $W\in\M_{k-s}$. By the same degree count $D\in\M_k$. With $R=UP_j+WP_{j+1}$\chg{, we have $D\beta-R=U\delta_j+W\delta_{j+1}$, where $|U\delta_j|<q^{s-\kappa}q^{-s}=q^{-\kappa}$ and $|W\delta_{j+1}|=q^{k-s-d_{j+2}}<q^{-\kappa}$, since $k+\kappa<s+d_{j+2}$. Hence $\nrm{D\beta}\le|D\beta-R|<q^{-\kappa}$. For uniqueness, note that $|D\beta-R|<1$ forces $R=[D\beta]$, so $R$ is determined by $D$; uniqueness of $(U,W)$ now follows from Lemma~\ref{lem:coord}.} The remaining statements follow by counting $U$ and $W$.
\end{proof}

\section{The Weil range: proof of Theorem~\ref{thm:weil}}\label{sec:weil}

In this section $q$ is an arbitrary prime power. It is convenient to work with the condition $\nrm{\alpha \chg{h}}\le q^{-M}$\chg{, where $h$ is a monic polynomial}.

\subsection{Reduction to residue classes}
\begin{lemma}[Reduction]\label{lem:reduction}
Let $a/f$ be in lowest terms with $f$ monic of degree $F$, and suppose $|\alpha-a/f|\le q^{-2F-g}$ with an integer $g\ge0$. Let $N\ge1$, $1\le M\le F$, $L=F-M$, and assume $N+M\le2F+g$. Then for every $\chg{h}\in\M_N$,
\[
\nrm{\alpha \chg{h}}\le q^{-M}\iff\deg(\overline{a\chg{h}})\le L,
\]
where $\overline{a\chg{h}}$ is the remainder of $a\chg{h}$ modulo $f$, and $\deg0=-\infty$.
\end{lemma}

\begin{proof}
Write $\alpha \chg{h}=a\chg{h}/f+\chg{h}(\alpha-a/f)$. Here $a\chg{h}/f=[a\chg{h}/f]+\overline{a\chg{h}}/f$, so $\nrm{a\chg{h}/f}=|\overline{a\chg{h}}|/|f|=q^{\deg\overline{a\chg{h}}-F}$. Also $|\chg{h}(\alpha-a/f)|\le q^{N-2F-g}\le q^{-M}$. By the ultrametric inequality\chg{, applied in both directions}, $\nrm{\alpha \chg{h}}\le q^{-M}$ iff $\nrm{a\chg{h}/f}\le q^{-M}$, i.e.\ iff $\deg\overline{a\chg{h}}\le F-M=L$.
\end{proof}

For $L\ge0$ put $W_L=\{b\in\A\setminus\{0\}:\deg b\le L\}$, $\Wp=W_L\cup\{0\}$ and $W_L^*(f)=\{b\in W_L:(b,f)=1\}$.

\begin{lemma}\label{lem:Wstar}
For every monic $f$ and $L\ge0$, $\bigl|\#W_L^*(f)-q^{L+1}\frac{\varphi(f)}{|f|}\bigr|\le2\cdot2^{\omega(f)}$.
\end{lemma}

\begin{proof}
For a divisor $d\mid f$, the set $\{b\in\Wp:d\mid b\}=\{dc:\deg c\le L-\deg d\}$ is an $\Fq$-vector space with $q^{L+1-\deg d}$ elements if $\deg d\le L$, and equals $\{0\}$ otherwise. By inclusion--exclusion over the squarefree divisors of $f$,
\[
\#W_L^*(f)=\sum_{\substack{d\mid f\\ \deg d\le L}}\mu(d)\bigl(q^{L+1-\deg d}-1\bigr)
=q^{L+1}\frac{\varphi(f)}{|f|}-q^{L+1}\sum_{\substack{d\mid f\\ \deg d>L}}\frac{\mu(d)}{|d|}-\sum_{\substack{d\mid f\\ \deg d\le L}}\mu(d),
\]
using $\sum_{d\mid f}\mu(d)/|d|=\varphi(f)/|f|$. Since $q^{L+1}/|d|\le1$ when $\deg d>L$, each of the last two terms is bounded by the number $2^{\omega(f)}$ of squarefree divisors of $f$.
\end{proof}

\subsection{The character method}
For $f$ monic of degree $F$, $(a,f)=1$, $N>F$ and $0\le L<F$, put
\[
S(f,a,N,L)=\#\{P\in\PP_N:\ \deg\overline{aP}\le L\}.
\]
Since $N>F$, every $P\in\PP_N$ is coprime to $f$, so $\overline{aP}$ is a unit modulo $f$ and the condition says $\overline{aP}\in W_L^*(f)$.

\begin{proposition}\label{prop:Sasymp}
In this situation,
\[
S(f,a,N,L)=\frac{\#W_L^*(f)}{\varphi(f)}\cdot\frac{q^N}{N}+O\Bigl(\frac{(F+1)\,q^{(N+L+1)/2}}{N}\Bigr).
\]
\end{proposition}

\begin{proof}
Distinct elements of $W_L$ are incongruent modulo $f$, since their degrees are less than $F$. By orthogonality of characters,
\[
S=\sum_{b\in W_L^*(f)}\#\{P\in\PP_N:\ aP\equiv b\bmod f\}=\sum_{b\in W_L^*(f)}\frac1{\varphi(f)}\sum_{\chi\bmod f}\overline\chi(b)\chi(a)\sum_{P\in\PP_N}\chi(P).
\]
The principal character contributes $\frac{\#W_L^*(f)}{\varphi(f)}\bigl(\frac{q^N}N+O(\frac{(F+1)q^{N/2}}N)\bigr)$ by Lemma~\ref{lem:PNT}, and $\#W_L^*(f)\le\varphi(f)$. By Lemma~\ref{lem:PNT}, the other characters contribute at most
\[
\frac{C(F+1)q^{N/2}}{N\varphi(f)}\sum_{\chi}\Bigl|\sum_{b\in W_L^*(f)}\overline\chi(b)\Bigr|
\le\frac{C(F+1)q^{N/2}}{N\varphi(f)}\,\varphi(f)^{1/2}\Bigl(\sum_\chi\Bigl|\sum_{b}\overline\chi(b)\Bigr|^2\Bigr)^{1/2},
\]
by Cauchy--Schwarz. By orthogonality the last sum equals $\varphi(f)\,\#W_L^*(f)\le\varphi(f)q^{L+1}$. This gives the error term.
\end{proof}

\subsection{The level principle}
\begin{definition}\label{def:admissible}
Let $0<\theta<1$ and $0<\lambda<1$. A sequence $(a_m/f_m)$ of fractions in lowest terms, with $f_m$ monic and $F_m=\deg f_m\to\infty$, is \emph{$(\theta,\lambda)$-admissible} if, for every $\varepsilon_0>0$, all large $m$ and all primes $N$ with $F_m\le\lambda N$ and $M:=\lceil\theta N\rceil\le(1-\varepsilon_0)F_m$, we have, with $L=F_m-M$,
\[
S(f_m,a_m,N,L)=\frac{\#W_L^*(f_m)}{\varphi(f_m)}\cdot\frac{q^N}{N}\bigl(1+o(1)\bigr)\chg{,}
\]
\chg{where $o(1)$ tends to $0$ as $m\to\infty$, uniformly in the primes $N$ in this range.}
\end{definition}

\begin{proposition}\label{prop:level}
Let $\nu\ge0$, $0<\theta<1$ and $0<\lambda<1$ with
\begin{equation}\label{eq:levelcond}
\theta<\lambda(2+\nu)-1\qquad\text{and}\qquad\theta(1+\nu)<1 .
\end{equation}
Suppose that $(a_m/f_m)$ is $(\theta,\lambda)$-admissible and $|\alpha-a_m/f_m|\le|f_m|^{-2-\nu}$ for all $m$. Then there are infinitely many monic irreducible $P$ with $\nrm{\alpha P}\le|P|^{-\theta}$.
\end{proposition}

\begin{proof}
Let $F=F_m$ be large and consider the interval $I_m=\bigl[\frac F\lambda,\ \frac{(2+\nu)F-1}{1+\theta}\bigr]$. The ratio of its endpoints tends to $\lambda(2+\nu)/(1+\theta)>1$ by \eqref{eq:levelcond}, so by the prime number theorem $I_m$ contains a prime $N$ for $m$ large. Put $M=\lceil\theta N\rceil$, $L=F-M$ and $g=\lfloor\nu F\rfloor$. We check the hypotheses of Lemma~\ref{lem:reduction} and of admissibility.
\begin{itemize}[leftmargin=2em]
\item $F\le\lambda N$, since $N\ge F/\lambda$.
\item $N+M\le(1+\theta)N+1\le(2+\nu)F$; as $N+M$ and $2F$ are integers, $N+M\le2F+g$.
\item $|\alpha-a/f|\le q^{-(2+\nu)F}\le q^{-2F-g}$.
\item $L\ge F-\frac{\theta((2+\nu)F-1)}{1+\theta}-1\ge\frac{1-\theta(1+\nu)}{1+\theta}F-1$, which tends to infinity by~\eqref{eq:levelcond}; and $L<F$ since $M\ge1$. Thus $M\le(1-\varepsilon_0)F$ with $\varepsilon_0=\frac{1-\theta(1+\nu)}{2(1+\theta)}$, for $m$ large.
\end{itemize}
Since $1\in W_L^*(f)$, admissibility gives $S(f,a,N,L)>0$ for $m$ large. By Lemma~\ref{lem:reduction}, every $P$ counted by $S$ satisfies $\nrm{\alpha P}\le q^{-M}\le q^{-\theta N}=|P|^{-\theta}$. Since $N\ge F/\lambda\to\infty$, we obtain infinitely many such $P$.
\end{proof}

\begin{corollary}\label{cor:charlevel}
For every $0<\theta<1$ and every $\lambda<1-\theta$, every sequence $(a_m/f_m)$ of fractions in lowest terms with $f_m$ monic and $\deg f_m\to\infty$ is $(\theta,\lambda)$-admissible.
\end{corollary}

\begin{proof}
Let $\varepsilon_0>0$, $F\le\lambda N$, $M=\lceil\theta N\rceil\le(1-\varepsilon_0)F$ and $L=F-M\ge\varepsilon_0F$. Note $N>F$. By Lemmas~\ref{lem:Wstar} and~\ref{lem:elem}(d),
\[
\frac{\#W_L^*(f)}{\varphi(f)}=\frac{q^{L+1}}{|f|}+O\Bigl(\frac{2^{\omega(f)}}{\varphi(f)}\Bigr)=q^{1-M}\bigl(1+O(q^{-L+o(F)})\bigr)=q^{1-M}(1+o(1)).
\]
Hence the main term in Proposition~\ref{prop:Sasymp} is $\asymp q^{N-M}/N$, and the ratio of the error term to the main term is
\[
\ll F\,q^{\frac{N+L+1}{2}-(N+1-M)}=F\,q^{\frac{F+M-N-1}{2}}\le F\,q^{\frac{(\lambda+\theta-1)N}{2}}\longrightarrow0,
\]
using $L+2M=F+M$, $F\le\lambda N$ and $M\le\theta N+1$.
\end{proof}

\subsection{Proof of Theorem~\ref{thm:weil}}
\begin{proof}[Proof of Theorem~\ref{thm:weil}]
(a) Let $\theta<\frac13$. Since $2(1-\theta)-1=1-2\theta>\theta$, we may choose $\lambda<1-\theta$ with $\theta<2\lambda-1$. The convergents of $\alpha$ satisfy $|\alpha-p_m/q_m|\le|q_m|^{-2}$ by Lemma~\ref{lem:cf}(c), and they form a $(\theta,\lambda)$-admissible sequence by Corollary~\ref{cor:charlevel}. Proposition~\ref{prop:level} with $\nu=0$ gives $\theta(\alpha)\ge\theta$.

(b) Let $0<\tau'<\tau=\min(\tau(\alpha),1)$. There are infinitely many $m$ with $\deg a_{m+1}\ge\tau'd_m$, and for these $|\alpha-p_m/q_m|=q^{-2d_m-\deg a_{m+1}}\le|q_m|^{-2-\tau'}$ by Lemma~\ref{lem:cf}(c). Let $\theta<\frac{1+\tau'}{3+\tau'}$. Then $(1-\theta)(2+\tau')-1>\theta$, so we may choose $\lambda<1-\theta$ with $\theta<\lambda(2+\tau')-1$. Moreover $\theta(1+\tau')<\frac{(1+\tau')^2}{3+\tau'}\le1$, since $(1+\tau')^2\le3+\tau'$ for $\tau'\le1$. By Corollary~\ref{cor:charlevel} the subsequence of these convergents is $(\theta,\lambda)$-admissible, and Proposition~\ref{prop:level} with $\nu=\tau'$ gives $\theta(\alpha)\ge\theta$. Let $\theta\to\frac{1+\tau'}{3+\tau'}$ and $\tau'\to\tau$.
\end{proof}

\begin{corollary}\label{cor:mahler}
Let $p\ge3$ be prime and $\alpha=\sum_{i\ge0}T^{-p^i}\in\mathbb F_p((1/T))$. Then $\alpha$ is algebraic of degree $p$ over $\Fq(T)$ for every power $q$ of $p$, and $\theta(\alpha)\ge\frac12$ over $\Fq$.
\end{corollary}

\begin{proof}
We have $\alpha^p=\sum_{i\ge1}T^{-p^i}=\alpha-T^{-1}$, so $\alpha$ is a root of $X^p-X+T^{-1}$. By Artin--Schreier theory, $X^p-X-c$ is irreducible over $K=\Fq(T)$ unless $c=y^p-y$ for some $y\in K$. For $c=-T^{-1}$ this is impossible: $c$ has a simple pole at $T=0$, whereas $y^p-y$ has a pole of order $pr$ at $T=0$ if $y$ has a pole of order $r\ge1$ there, and no pole otherwise. So $\alpha$ has degree $p$, and in particular $\alpha\notin\Fq(T)$.

For $k\ge0$, the partial sum $r_k=\sum_{i\le k}T^{-p^i}$ equals $u_k/T^{p^k}$ with $u_k=\sum_{i\le k}T^{p^k-p^i}$. The term $i=k$ of $u_k$ is $1$, so $u_k$ is coprime to $T$. Moreover $|\alpha-r_k|=|T|^{-p^{k+1}}=|T^{p^k}|^{-p}<|T^{p^k}|^{-2}$. By Lemma~\ref{lem:legendre}, $r_k$ is a convergent, and by Lemma~\ref{lem:cf}(c) the following partial quotient has degree $(p-2)p^k$. Hence $\tau(\alpha)\ge p-2\ge1$, and Theorem~\ref{thm:weil}(b) gives $\theta(\alpha)\ge\frac12$.
\end{proof}

\begin{remark}\label{rem:lacunary}
More generally, Theorem~\ref{thm:weil}(b) gives $\theta(\alpha)\ge\frac12$ for every $\alpha$ with $\tau(\alpha)\ge1$, over every $\Fq$, including prime $q$. For example, let $\alpha=\sum_{i\ge0}c_iT^{-n_i}$ with $c_i\in\Fq^\times$ and $n_{i+1}\ge3n_i$ for infinitely many $i$, and $n_{i+1}>2n_i$ for all $i$. The argument above shows that the partial sums $\sum_{i\le k}c_iT^{-n_i}$ are convergents of $\alpha$, with following partial quotient of degree $n_{k+1}-2n_k$, so that $\tau(\alpha)\ge1$. Such $\alpha$ is not in $\Fq(T)$, since the coefficients of an element of $\Fq(T)$ in $\Fq((1/T))$ are eventually periodic.
\end{remark}

\section{Squarefree moduli: proof of Theorem~\ref{thm:sqfree}}\label{sec:sqfree}

In this section $q$ is again an arbitrary prime power. We first recall the theorem of Sawin which we shall use, then deduce from it a level of distribution for primes in the Diophantine sets, and finally show that the squarefree hypothesis of Theorem~\ref{thm:sqfree} is satisfied for most $\alpha$.

\begin{theorem}[Sawin {\cite[Theorem~1.2]{Sawin}}]\label{thm:sawin}
Let $0<\lambda<1$ and let $q$ be a prime power with $q>e^2\,4^{\frac1{1-\lambda}}\bigl(\frac1{1-\lambda}-\frac12\bigr)^2$. There is $\delta=\delta(q,\lambda)>0$ such that, for all integers $N,F$ with $F\le\lambda N$, every squarefree $g\in\M_F$ and every $c\in(\A/g)^\times$,
\[
\sum_{\substack{\chg{h}\in\M_N\\ \chg{h}\equiv c\bmod g}}\Lambda(\chg{h})=\frac{q^N}{\varphi(g)}+O_{q,\lambda}\bigl(q^{(N-F)(1-\delta)}\bigr).
\]
\end{theorem}

The condition on $q$ is exactly $\lambda<\lambda(q)$, by~\eqref{eq:lambdaq} and the monotonicity of the right side of~\eqref{eq:lambdaq} in $\lambda$. For instance, $\lambda=\frac12$ gives $q>36e^2=266.0\ldots$, which is the threshold for a level of distribution beyond $\frac12$ stated in~\cite[\S1]{Sawin}.

\begin{corollary}\label{cor:sawinlevel}
Let $q$ and $\lambda$ be as in Theorem~\ref{thm:sawin}, and $0<\theta<1$. Every sequence $(a_m/f_m)$ of fractions in lowest terms with $f_m$ monic and squarefree and $\deg f_m\to\infty$ is $(\theta,\lambda)$-admissible \textup{(}Definition~\ref{def:admissible}\textup{)}.
\end{corollary}

\begin{proof}
Let $f=f_m$, $a=a_m$, $F=\deg f$, and let $N$ be prime with $F\le\lambda N$, $M=\lceil\theta N\rceil\le(1-\varepsilon_0)F$ and $L=F-M$. Since $N>F$, every $P\in\PP_N$ is coprime to $f$, and $P$ is counted by $S(f,a,N,L)$ iff $P\equiv\bar a b\bmod f$ for some $b\in W_L^*(f)$, where $\bar a$ is the inverse of $a$ modulo $f$. Since $N$ is prime, the only prime powers $P^e\in\M_N$ with $e\ge2$ are $P^N$ with $\deg P=1$; there are $q$ of them, each with $\Lambda(P^N)=1$. Hence, by Theorem~\ref{thm:sawin},
\[
N\cdot S(f,a,N,L)=\sum_{b\in W_L^*(f)}\ \sum_{\substack{\chg{h}\in\M_N\\ \chg{h}\equiv\bar ab\bmod f}}\Lambda(\chg{h})+O(q)
=\#W_L^*(f)\Bigl(\frac{q^N}{\varphi(f)}+O\bigl(q^{(N-F)(1-\delta)}\bigr)\Bigr)+O(q).
\]
Since $q^N/\varphi(f)\ge q^{N-F}$, the relative error of the first term is $O(q^{-\delta(N-F)})$, which tends to $0$ because $N-F\ge(1-\lambda)N$. By Lemma~\ref{lem:Wstar}, $\#W_L^*(f)\,q^N/\varphi(f)\ge q^{N-M}(1+o(1))\to\infty$, so the term $O(q)$ is also negligible.
\end{proof}

\begin{lemma}\label{lem:pencil}
Let $m\ge1$ and $c\in\Fq$, and put $a=p_m+cp_{m-1}$, $f=q_m+cq_{m-1}$. Then $f$ is monic of degree $d_m$, $(a,f)=1$ and $|\alpha-a/f|\le|f|^{-2}$.
\end{lemma}

\begin{proof}
Since $d_{m-1}<d_m$, $f$ is monic of degree $d_m$. By Lemma~\ref{lem:cf}(a) (for the normalised convergents, \eqref{eq:cf}), $af'-fa'$ with $(a',f')=(p_{m-1},q_{m-1})$ equals $p_mq_{m-1}-q_mp_{m-1}\in\Fq^\times$, so $(a,f)=1$. Finally $f\alpha-a=\delta_m+c\delta_{m-1}$ has absolute value at most $\max(q^{-d_{m+1}},q^{-d_m})=|f|^{-1}$, so $|\alpha-a/f|\le|f|^{-2}$.
\end{proof}

\begin{proposition}\label{prop:measure}
Normalise the Haar measure on $\mathfrak m=\{x\in\Kinf:|x|<1\}$ to total mass $1$. For $D\ge2$, let $E_D\subset\mathfrak m$ be the set of $\alpha$ having a squarefree convergent denominator of degree $D$. Then $\operatorname{vol}(E_D)\ge1-\frac2q$. Consequently, the set of $\alpha\in\mathfrak m$ with infinitely many squarefree convergent denominators has measure at least $1-\frac2q$.
\end{proposition}

\begin{proof}
Let $r/Q$ be in lowest terms with $Q\in\M_D$. If $|\alpha-r/Q|<q^{-2D}$, then $r/Q$ is a convergent of $\alpha$ by Lemma~\ref{lem:legendre}. Conversely, if $r/Q$ is a convergent, then $|\alpha-r/Q|=q^{-d_i-d_{i+1}}\le q^{-2D-1}$ by Lemma~\ref{lem:cf}(c). So $E_D$ is the disjoint union of the balls $\{|\alpha-r/Q|<q^{-2D}\}$ over $Q\in\M_D$ squarefree and $r$ coprime to $Q$ with $\deg r<D$; such a ball lies in $\mathfrak m$ and has measure $q^{-2D}$. The balls are disjoint, because $|r/Q-r'/Q'|\ge|QQ'|^{-1}=q^{-2D}$ for distinct fractions. Hence
\[
\operatorname{vol}(E_D)=q^{-2D}\sum_{\substack{Q\in\M_D\\ Q\text{ squarefree}}}\varphi(Q)\ \ge\ q^{-2D}\Bigl(\sum_{Q\in\M_D}\varphi(Q)-\sum_{\substack{Q\in\M_D\\ Q\text{ not squarefree}}}|Q|\Bigr).
\]
From the generating functions
\[
\sum_Q\varphi(Q)u^{\deg Q}=\frac{1-qu}{1-q^2u},\qquad\sum_Q\mu^2(Q)u^{\deg Q}=\frac{1-qu^2}{1-qu}
\]
we get $\sum_{Q\in\M_D}\varphi(Q)=q^{2D}(1-q^{-1})$, and there are exactly $q^{D-1}$ non-squarefree $Q\in\M_D$ for $D\ge2$. Hence $\operatorname{vol}(E_D)\ge1-\frac1q-\frac1q$. The last statement follows from $\operatorname{vol}(\limsup E_D)\ge\limsup\operatorname{vol}(E_D)$, which is the reverse Fatou lemma for the probability measure $\operatorname{vol}$.
\end{proof}

\begin{proof}[Proof of Theorem~\ref{thm:sqfree}]
Let $\lambda<\lambda(q)$ and $\theta<2\lambda-1$. For each of the infinitely many $m$ in the hypothesis, Lemma~\ref{lem:pencil} gives a fraction $a/f$ in lowest terms, with $f$ monic, squarefree, $\deg f=d_m\to\infty$, and $|\alpha-a/f|\le|f|^{-2}$. By Corollary~\ref{cor:sawinlevel} these fractions form a $(\theta,\lambda)$-admissible sequence, and Proposition~\ref{prop:level} with $\nu=0$ gives $\theta(\alpha)\ge\theta$. Letting $\lambda\to\lambda(q)$ gives $\theta(\alpha)\ge2\lambda(q)-1$. The measure statement is Proposition~\ref{prop:measure} with $c=0$; it applies in every translate of $\mathfrak m$, since translation by $\A$ does not change the convergent denominators.
\end{proof}

\begin{remark}\label{rem:sqfree}
(a) Theorem~\ref{thm:sqfree} uses Proposition~\ref{prop:level} only with $\nu=0$. If $\nu>0$ and there are infinitely many $m$ such that $q_m$ is squarefree and $\deg a_{m+1}\ge\nu d_m$, then $|\alpha-p_m/q_m|\le|q_m|^{-2-\nu}$ by Lemma~\ref{lem:cf}(c), and Proposition~\ref{prop:level} together with Corollary~\ref{cor:sawinlevel} gives
\[
\theta(\alpha)\ge\min\Bigl(\lambda(q)(2+\nu)-1,\ \frac1{1+\nu}\Bigr).
\]

(b) The bound $1-\frac2q$ in Proposition~\ref{prop:measure} is certainly not optimal. The events ``$\alpha$ has a squarefree convergent denominator of degree $D$'' should be nearly independent for well-separated $D$, and one expects, via a quasi-independent Borel--Cantelli argument, that almost every $\alpha$ has infinitely many squarefree convergent denominators. We do not pursue this here.
\end{remark}

\section{\texorpdfstring{The M\"obius function along linear sequences}{The Moebius function along linear sequences}}\label{sec:SS}

This section contains the input on the M\"obius function which drives Theorems~\ref{thm:main} and~\ref{thm:main2}. We quote two results of Sawin and Shusterman, and then apply them to the fibres of Lemma~\ref{lem:fibre}. Both rest on Pellet's formula, which expresses $\mu(f)$ through the discriminant of $f$, and on the fact that in characteristic $p$ the polynomials $r+s^p$ all have the same derivative. It follows that $\mu(r+s^p)$ is essentially a quadratic Dirichlet character evaluated at a shift of $s$, and short sums of such characters are estimated by geometric methods.

To avoid clashes with the notation of the rest of the paper, we write the parameters of these results as $\epsilon_1,\epsilon_2,\psi$ (Theorem~\ref{thm:SS}) and $\zeta,\Upsilon_\zeta,\varpi,\ell$ (Theorem~\ref{thm:SS2}); in~\cite{SS,SS2} they are called $\epsilon,\delta,\beta$ and $\gamma,\beta,\alpha,c_1$, respectively (with $k=1$ in~\cite[Theorem~5.1]{SS2}, so that $\beta=(1+2\gamma)^k=1+2\gamma$).

\begin{theorem}[Sawin--Shusterman {\cite[Theorem~4.5]{SS}}]\label{thm:SS}
Let $p$ be an odd prime, $\epsilon_1,\epsilon_2>0$, $0<\psi<\frac12$, and let $q$ be a power of $p$ with
\[
q>\Bigl(\frac{pe}{\eta_0}\Bigr)^{\frac2{1-2\psi}},\qquad\eta_0=\min\Bigl\{\frac{\epsilon_1}{\epsilon_1+2},\ \frac{\epsilon_1\epsilon_2}{\epsilon_1+\epsilon_2}\Bigr\}.
\]
Let $t,m,k'\ge0$ be integers with $t\ge\max(\epsilon_1m,\epsilon_2k')$ and $k'\ne t+m$, and let $a\in\M_{k'}$, $M_0\in\M_m$. Then
\[
\sum_{g\in\M_t}\mu(a+gM_0)\ll q^{t(1-\psi/p)}\qquad(t\to\infty),
\]
where the implied constant depends only on $\epsilon_1,\epsilon_2,\psi,q$.
\end{theorem}

\begin{remark*}
In its general form,~\cite[Theorem~4.5]{SS} concerns sums $\sum_{g\in\M_d}\prod_{i=1}^{n}\mu(a_i+gM_i)$ with $(a_i,M_i)\in\M_{k_i}\times\M_{m_i}$, $d\ge\max(\epsilon m_i,\delta k_i)$, $k_i\ne d+m_i$ and pairwise distinct fractions $a_i/M_i$, under the condition
\[
q>\Bigl(\frac{pne}{\min\{\frac{\epsilon}{\epsilon+2},\frac{\epsilon\delta}{\epsilon+\delta}\}}\Bigr)^{\frac2{1-2\beta}}
\]
of~\cite[(4.33)]{SS}. It bounds them by $O(|\M_d|^{1-\beta/p})=O(q^{d(1-\beta/p)})$ as $d\to\infty$, with an implied constant depending only on $\epsilon,\delta,\beta,n$ and $q$. Theorem~\ref{thm:SS} is the case $n=1$, where $pne=pe$ and the distinctness hypothesis is vacuous. We shall use the following three features.
\begin{enumerate}[label=\rm(\roman*)]
\item The implied constant does not depend on $a$ or $M_0$. This uniformity is essential: in Proposition~\ref{prop:goodfibre} both $a$ and $M_0$ depend on $B$ through the convergents of $B\alpha$, and we sum over all $B$.
\item The condition $k'\ne t+m$ excludes cancellation between the leading terms of $a$ and $gM_0$. In our application $k'=m<t+m$, so every $a+gM_0$ is monic of degree $t+m$.
\item Neither $M_0$ squarefree nor $(a,M_0)=1$ is assumed. At the end of the proof of~\cite[Proposition~4.3]{SS}, the case $H=\gcd(a,M_0)\ne1$ is reduced to the coprime one. If $H$ is not squarefree the sum vanishes; otherwise one factors out $\mu(H)$ and removes residue classes at a cost $q^{o(t)}$ \chg{(see~\cite[(4.31)--(4.32)]{SS})}, which is absorbed by slightly increasing $\psi$. This is possible because the condition on $q$ is a strict inequality.
\end{enumerate}
\chg{As a consequence,}~\cite[Theorem~1.7]{SS} gives a level of distribution of $\mu$ close to $1$ in arbitrary progressions; compare~\cite[\S1]{Sawin}, where the analogue for squarefree moduli is proved with a condition on $q$ independent of $p$.
\end{remark*}

\begin{theorem}[Sawin--Shusterman {\cite[Theorem~5.1]{SS2}}]\label{thm:SS2}
Let $p$ be an odd prime and $q$ a power of $p$. Let $0<\zeta\le1$, put $\Upsilon_\zeta=1+2\zeta$, and let
\[
0<\varpi<\frac1{2p}+\frac{\log_q\zeta}{p}-\log_q(1+2\zeta).
\]
Let \chg{$\ell\ge0$ be an integer and let} $a_0,a_1\in\A$ with $a_1\ne0$ and $\deg a_0,\deg a_1\le\ell$, and let $\I=\{f_0+g:\ \deg g<d\}$ with $\deg f_0\le d$. Then
\[
\sum_{f\in\I}\mu(a_0+a_1f)\ll q^{d(1-\varpi)}\,\Upsilon_\zeta^{2\ell}\qquad(d\to\infty),
\]
with the implied constant depending only on $q,\varpi,\zeta$.
\end{theorem}

\begin{proof}[Deduction]
We recall~\cite[Theorem~5.1]{SS2}, in the notation used there: polynomials are in a variable $u$, which we identify with our $T$, and $T$ is a formal variable. Let $k\ge1$, $0<\gamma\le1$ and $0<\alpha<\frac1{2p}+\frac{\log_q\gamma}p-k\log_q(1+2\gamma)$, and put $\beta=(1+2\gamma)^k$. Let $F(u,T)=\sum_{i=0}^ka_i(u)T^i\in\Fq[u,T]$ be separable of degree $k$ in $T$, and let $c_1\ge0\ge c_2$ with $\deg a_i\le c_1+c_2i$ for $0\le i\le k$. Let $g\in\Fq[u]$ be squarefree, $t$ an infinitame $g$-periodic trace function, and $\I$ an interval as in~\cite[Notation~4.1]{SS2}. Then
\begin{multline*}
\sum_{f\in\I}\mu(F(u,f))\,t(f)\ll q^{\dim(\I)(1-\alpha)}\beta^{2c_1+(k+1)c_2}\bigl(\beta^{-c_2-\dim(\I)}+\beta^{\deg(\I)-\dim(\I)}\bigr)\\
\times\prod_{\pi\mid g}\bigl(r(t_\pi)(1+\gamma)+c(t_\pi)\gamma\bigr)^{\deg\pi}
\end{multline*}
as $\dim(\I)\to\infty$, with the implied constant depending only on $q,k,\alpha,\gamma$.

We apply this with $k=1$, $F(u,T)=a_0+a_1T$ (separable, of degree $1$ in $T$, since $a_1\ne0$), $c_1=\ell$, $c_2=0$, the modulus $g=1$ and the trace function $t\equiv1$ \chg{(the trace function of the constant sheaf, which is infinitame)}; then $\beta=1+2\zeta=\Upsilon_\zeta$, the condition on $\alpha=\varpi$ is the one in Theorem~\ref{thm:SS2}, $\deg a_0,\deg a_1\le\ell=c_1+c_2i$, and the product over $\pi\mid g$ is empty. In the notation of~\cite[Notation~4.1]{SS2}, $q^{\dim\I}=\#\I=q^d$ and \chg{$\deg\I\le d=\dim\I$, since every element of $\I$ has degree at most $d$; hence $\Upsilon_\zeta^{\deg\I-\dim\I}\le1$, and} the bound reads $q^{d(1-\varpi)}\Upsilon_\zeta^{2\ell}(\Upsilon_\zeta^{-d}+\chg{\Upsilon_\zeta^{\deg\I-\dim\I}})\le2q^{d(1-\varpi)}\Upsilon_\zeta^{2\ell}$, with an implied constant depending only on $q,\varpi,\zeta$.
\end{proof}

The bound of Theorem~\ref{thm:SS2} is uniform in $a_0,a_1$ and $f_0$, and depends on the coefficients only through the factor $\Upsilon_\zeta^{2\ell}$, which grows with the degree bound $\ell$. This factor causes no loss as long as $\ell$ is at most a fixed multiple of the length $d$ of the interval. In Proposition~\ref{prop:twodir} we take $\ell=k$ and $d\ge\eta k$, so that $\Upsilon_\zeta^{2\ell}\le q^{(2/\eta)d\log_q\Upsilon_\zeta}$; this is why the margin~\eqref{eq:margin} below involves $\frac2\eta\log_q\Upsilon_{\zeta_0}$.

\begin{proposition}\label{prop:goodfibre}
Let $p,\varepsilon,q$ be as in Theorem~\ref{thm:main}. Fix $\psi=\psi(q,\varepsilon)\in(0,\frac12)$ such that the hypothesis of Theorem~\ref{thm:SS} holds with $\epsilon_1=\epsilon_2=\varepsilon$; this is possible since $\eta_0=\frac{\varepsilon}{\varepsilon+2}$ and~\eqref{eq:qcond} says $q>(pe/\eta_0)^2$ with strict inequality. There is $C=C(q,\varepsilon)$ with the following property. Let $\gamma\in\Kinf\setminus\Fq(T)$, $k\ge1$, $\kappa\ge0$, and let $j=j_\gamma(k+\kappa)\ge0$ and $s=d_{j+1}$ satisfy $\kappa\le s\le\frac{k}{1+\varepsilon}$. Then, with $\F=\{D\in\M_k:\ \nrm{D\gamma}<q^{-\kappa}\}$ and $t=k-s$,
\[
\Bigl|\sum_{D\in\F}\mu(D)\Bigr|\le C\,q^{\,k-\kappa-\psi t/p}.
\]
\end{proposition}

\begin{proof}
By Lemma~\ref{lem:fibre}, $\sum_{D\in\F}\mu(D)=\sum_{|U|<q^{s-\kappa}}\Sigma_U$, with $\Sigma_U=\sum_{W\in\M_t}\mu(UQ_j+WQ_{j+1})$, and there are $q^{s-\kappa}$ choices of $U$. We have $t\ge k-\frac{k}{1+\varepsilon}=\frac{\varepsilon k}{1+\varepsilon}\ge\varepsilon s$; in particular $t\ge1$ once $k\ge3/\varepsilon$, and for smaller $k$ the claim holds \chg{with $C=q^{3/\varepsilon}$, by the trivial bound $|\F|\le q^{k-\kappa}$ and $\psi t/p<k<3/\varepsilon$}.

Fix $U$ and divide: $UQ_j=hQ_{j+1}+a_0$ with $\deg a_0<s$. If $UQ_j\ne0$, then $\deg h=\deg(UQ_j)-s<(s-\kappa+d_j)-s=d_j-\kappa\le k-s=t$, by~\eqref{eq:jdef}; otherwise $h=0$. Put $a=a_0+Q_{j+1}$, which is monic of degree $s$, and $g=W+h-1$. Since $\deg(h-1)<t$, the map $W\mapsto g$ is a bijection of $\M_t$, and
\[
a+gQ_{j+1}=a_0+(W+h)Q_{j+1}=UQ_j+WQ_{j+1}.
\]
Thus $\Sigma_U=\sum_{g\in\M_t}\mu(a+gQ_{j+1})$. We apply Theorem~\ref{thm:SS} with $M_0=Q_{j+1}$, $m=k'=s$: the conditions $t\ge\max(\epsilon_1m,\epsilon_2k')=\varepsilon s$ and $k'\ne t+m$ hold. This gives $|\Sigma_U|\le C\,q^{t(1-\psi/p)}$, with $C$ independent of $U,\gamma,k,\kappa$ by the uniformity in $a$ and $M_0$ (see the remark after Theorem~\ref{thm:SS}); enlarging $C$ covers the finitely many small $t$. Summing over $U$ gives $C\,q^{s-\kappa+t-\psi t/p}=C\,q^{k-\kappa-\psi t/p}$.
\end{proof}

\begin{proposition}\label{prop:twodir}
Let $p$ be an odd prime, $0<\eta\le\frac1{12}$, and $q$ a power of $p$ satisfying~\eqref{eq:qcond2}. There are $C=C(q,\eta)$ and $c_0=c_0(q,\eta)>0$ with the following property. Let $\gamma\in\Kinf\setminus\Fq(T)$, $k\ge1$, $\kappa\ge0$, and let $j=j_\gamma(k+\kappa)\ge0$ and $s=d_{j+1}$ satisfy $\kappa\le s\le k$ and
\[
\max(k-s,\ s-\kappa)\ge\eta k .
\]
Then, with $\F=\{D\in\M_k:\ \nrm{D\gamma}<q^{-\kappa}\}$,
\[
\Bigl|\sum_{D\in\F}\mu(D)\Bigr|\le C\,q^{\,k-\kappa-c_0\eta k}.
\]
\end{proposition}

\begin{proof}
\emph{Step 1: choice of parameters.} Put $\Gamma=2p(1+2/\eta)$ and $\zeta_0=1/\Gamma\le1$, and let
\[
\varpi_{\max}=\frac1{2p}+\frac{\log_q\zeta_0}{p}-\log_q(1+2\zeta_0),\qquad
\Delta=\varpi_{\max}-\frac2\eta\log_q(1+2\zeta_0).
\]
Since $\log(1+2\zeta_0)\le2\zeta_0$, we have $(1+\frac2\eta)\log_q(1+2\zeta_0)\le\frac{2(1+2/\eta)}{\Gamma\log q}=\frac{1}{p\log q}$, and hence
\[
\Delta=\frac1{2p}-\frac{\log\Gamma}{p\log q}-\Bigl(1+\frac2\eta\Bigr)\log_q(1+2\zeta_0)\ \ge\ \frac{1}{p\log q}\Bigl(\tfrac12\log q-\log\Gamma-1\Bigr).
\]
This is positive because~\eqref{eq:qcond2} says $q>e^2\Gamma^2$. Put $c_0=\Delta/2$ and $\varpi=\varpi_{\max}-c_0$. Then $0<\varpi<\varpi_{\max}$ (as $\varpi=c_0+\frac2\eta\log_q(1+2\zeta_0)>0$), so $\varpi$ satisfies the hypothesis of Theorem~\ref{thm:SS2} with $\zeta=\zeta_0$, and
\begin{equation}\label{eq:margin}
\varpi-\frac2\eta\log_q\Upsilon_{\zeta_0}=\Delta-c_0=c_0 .
\end{equation}

\emph{Step 2: structure.} By Lemma~\ref{lem:fibre}, $\F=\{UQ_j+WQ_{j+1}:|U|<q^{s-\kappa},\ W\in\M_t\}$ with $t=k-s$, and every element is monic of degree $k$. Let $C'$ be the implied constant of Theorem~\ref{thm:SS2} for $(q,\varpi,\zeta_0)$, valid for $d\ge d_0$.

\emph{Step 3: case $t\ge\eta k$.} For fixed $U$, apply Theorem~\ref{thm:SS2} with $a_0=UQ_j$ (of degree $<k$), $a_1=Q_{j+1}$ (of degree $s\le k$), $\ell=k$, and $\I=\M_t=\{T^t+g:\deg g<t\}$, so $d=t$. Since $k/t\le1/\eta$, \eqref{eq:margin} gives
\[
\Bigl|\sum_{W\in\M_t}\mu(UQ_j+WQ_{j+1})\Bigr|\le C'q^{t(1-\varpi)+2k\log_q\Upsilon_{\zeta_0}}=C'q^{t-(\varpi-\frac{2k}{t}\log_q\Upsilon_{\zeta_0})t}\le C'q^{t-c_0t}.
\]
Summing over the $q^{s-\kappa}$ values of $U$ gives $C'q^{k-\kappa-c_0t}\le C'q^{k-\kappa-c_0\eta k}$.

\emph{Step 4: case $s-\kappa\ge\eta k$.} For fixed $W$, apply Theorem~\ref{thm:SS2} with $a_0=WQ_{j+1}$ (of degree $k$), $a_1=Q_j\ne0$ (of degree $d_j\le k$), $\ell=k$, and $\I=\{U:\deg U<s-\kappa\}$ (so $f_0=0$ and $d=s-\kappa$). Since $k/d\le1/\eta$, the inner sum is at most $C'q^{(s-\kappa)(1-c_0)}$ as in Step~3. Summing over the $q^t$ values of $W$ gives $C'q^{k-\kappa-c_0(s-\kappa)}\le C'q^{k-\kappa-c_0\eta k}$.

\emph{Step 5.} In both cases $d\ge\eta k$. If $\eta k<d_0$, the trivial bound $|\F|=q^{k-\kappa}$ suffices after enlarging $C$. So the claim holds with $C=\max(C',q^{c_0d_0})$.
\end{proof}

\section{Proof of Theorem~\ref{thm:main}}\label{sec:proof}
We now prove Theorem~\ref{thm:main}, following the plan described in the introduction. The small divisors are handled exactly by Lemma~\ref{lem:hankel}, and the large divisors by the fibre decomposition of Lemma~\ref{lem:fibre} combined with Proposition~\ref{prop:goodfibre}. It may help to keep in mind the case of bounded partial quotients, in which $E_\alpha(X)=O(1)$ and condition~\eqref{eq:range} simply reads $\kappa\le(\frac12-3\varepsilon)n+O(1)$.

Fix $\alpha$, $n$ and $\kappa$ satisfying~\eqref{eq:range}, and put $E=E_\alpha(n+\kappa)$ and $\V=\V_n(\alpha,\kappa)$. Every application of Lemma~\ref{lem:count} below, namely in \S\ref{sec:sizeV}, in Proposition~\ref{prop:small} and in case~(b) of Proposition~\ref{prop:large}, has $L+M=n+\kappa$ exactly, so it may be made with this $E$. Condition~\eqref{eq:range} reads
\begin{equation}\label{eq:rangeE}
\kappa+\tfrac E2\le\bigl(\tfrac12-3\varepsilon\bigr)n,
\end{equation}
and it implies $E\le n-2\kappa-6\varepsilon n$ and $\kappa\le(\frac12-3\varepsilon)n$.

\subsection{\texorpdfstring{The size of $\V$ and the decomposition}{The size of V and the decomposition}}\label{sec:sizeV}
By Lemma~\ref{lem:hankel} with $\beta=\alpha$ and $N=n$, we have $|\V|=q^{n-\kappa}$ unless there is $G\ne0$ with $\deg G<\kappa$ and $\nrm{G\alpha}<q^{-n}$. Lemma~\ref{lem:count}(ii), applied with $L=\kappa<M=n$, shows that there is no such $G$, since $M-L=n-\kappa\ge E$. Hence $|\V|=q^{n-\kappa}$.

By Lemma~\ref{lem:elem}(c),
\begin{equation}\label{eq:decomp}
S:=\sum_{f\in\V}\Lambda(f)=-\sum_{f\in\V}\sum_{D\mid f}\mu(D)\deg D=-\sum_{k=1}^nk\,\Sigma_k,\qquad\Sigma_k=\sum_{D\in\M_k}\mu(D)\,C(D),
\end{equation}
where $C(D)=\#\{B\in\M_{n-k}:\ DB\in\V\}=\#\{B\in\M_{n-k}:\ \nrm{B\cdot D\alpha}<q^{-\kappa}\}$. (The term $k=0$ has weight $0$.)

\subsection{Small divisors}
\begin{proposition}\label{prop:small}
Let $E'\ge E_\alpha(n+\kappa)$. For $1\le k\le n-\kappa$,
\[
\sum_{D\in\M_k}\bigl|C(D)-q^{n-k-\kappa}\bigr|\le q^{o(n)}\bigl(q^{k+E'/2}+q^{(n-\kappa+E')/2}\bigr).
\]
\end{proposition}

\begin{proof}
Fix $D\in\M_k$, and apply Lemma~\ref{lem:hankel} with $\beta=D\alpha$ and $N=n-k\ge\kappa$. It gives
\[
\bigl|C(D)-q^{n-k-\kappa}\bigr|\le q^{n-k-\kappa}\,\#\{G\ne0:\ \deg G<\kappa,\ \nrm{GD\alpha}<q^{-(n-k)}\}.
\]
Summing over $D$, we see that the left side of the claim is at most $q^{n-k-\kappa}$ times the number of pairs $(D,G)$ with $D\in\M_k$, $G\ne0$, $\deg G<\kappa$ and $\nrm{GD\alpha}<q^{-(n-k)}$.

We count these pairs through the product $g=GD$. Then $g\ne0$, $\deg g<\kappa+k\le n$ and $\nrm{g\alpha}<q^{-(n-k)}$. Conversely, $g$ determines the pair once the monic divisor $D$ of $g$ is chosen, so each $g$ arises from at most $\tau_d(g)=q^{o(n)}$ pairs. Hence the number of pairs is at most $q^{o(n)}\mathcal N(L,M)$, with $L=\kappa+k$ and $M=n-k$, so that $L+M=n+\kappa$. We now apply Lemma~\ref{lem:count}. If $L\ge M$, then $\mathcal N\le q^{2k+\kappa-n+E'/2}$, and if $L<M$, then $\mathcal N\le q^{(E'+2k+\kappa-n)/2}$. Multiplying by $q^{n-k-\kappa}$, we obtain the two bounds
\begin{align*}
q^{n-k-\kappa}\,q^{2k+\kappa-n+E'/2+o(n)}&=q^{k+E'/2+o(n)}&&(L\ge M),\\
q^{n-k-\kappa}\,q^{(E'+2k+\kappa-n)/2+o(n)}&=q^{(n-\kappa+E')/2+o(n)}&&(L<M),
\end{align*}
which together give the proposition.
\end{proof}

By Lemma~\ref{lem:elem}(b), $\sum_{D\in\M_k}\mu(D)q^{n-k-\kappa}$ equals $-q^{n-\kappa}$ for $k=1$ and $0$ for $k\ge2$. Hence, for every integer $1\le K\le n-\kappa$ and every $E'\ge E_\alpha(n+\kappa)$,
\begin{equation}\label{eq:small}
-\sum_{k\le K}k\,\Sigma_k=q^{n-\kappa}+O\Bigl(n^2q^{o(n)}\bigl(q^{K+E'/2}+q^{(n-\kappa+E')/2}\bigr)\Bigr).
\end{equation}
In this section we take $K=\lfloor n/2\rfloor\le n-\kappa$ and $E'=E$.

\subsection{Large divisors}
\begin{proposition}\label{prop:large}
For $n/2<k\le n$,
\[
|\Sigma_k|\ll q^{\,n-\kappa-\psi\varepsilon n/(3p)}+q^{o(n)}\bigl(q^{n/2+2\varepsilon n+E/2}+q^{(n-\kappa+E)/2+\varepsilon n}\bigr),
\]
with an implied constant depending only on $q$ and $\varepsilon$.
\end{proposition}

\begin{proof}
Interchanging the summations, we have
\[
\Sigma_k=\sum_{B\in\M_{n-k}}\ \sum_{D\in\F_B}\mu(D),\qquad\F_B=\{D\in\M_k:\ \nrm{D\cdot B\alpha}<q^{-\kappa}\}.
\]
For each $B$ let $\gamma=B\alpha$, which is not in $\Fq(T)$. We use its convergent denominators $Q_i=Q_i(B)$, of degrees $d_i=d_i(B)$, and we put $j=j_\gamma(k+\kappa)$ and $s=d_{j+1}$. Note that $k>n/2\ge\kappa$. In fact $\frac{k}{1+\varepsilon}>\frac n{2+2\varepsilon}\ge\kappa$, since $(\frac12-3\varepsilon)(2+2\varepsilon)\le1$. We split the $B$ into four classes, according to the values of $j$ and $s$.

\emph{(o) $j=-1$.} Then $\F_B=\M_k$ by Lemma~\ref{lem:fibre}(a), and $\sum_{D\in\M_k}\mu(D)=0$ since $k\ge2$ (recall $k>n/2$ and $n\ge n_0\ge3$).

\emph{(a) $j\ge0$ and $\kappa\le s\le k/(1+\varepsilon)$.} Then $t=k-s\ge\frac{\varepsilon k}{1+\varepsilon}\ge\frac{\varepsilon n}{2(1+\varepsilon)}\ge\frac{\varepsilon n}3$, and Proposition~\ref{prop:goodfibre} bounds the inner sum by $Cq^{k-\kappa-\psi\varepsilon n/(3p)}$. Summing over the $q^{n-k}$ values of $B$ gives the first term of the claim.

\emph{(b) $j\ge0$ and $s>k/(1+\varepsilon)$.} If $s>k$, then $\F_B=\emptyset$. Otherwise $\kappa<s\le k$, and $|\F_B|=q^{k-\kappa}$ by Lemma~\ref{lem:fibre}. Here we do not look for cancellation; instead we count these $B$. Fix the value of $s$; there are at most $n+\kappa+1$ possible values. Put $g=Q_jB$, which is monic and nonzero. By~\eqref{eq:jdef} and~\eqref{eq:cf},
\[
\deg g=d_j+n-k\le(k+\kappa-s)+n-k=n+\kappa-s,\qquad\nrm{g\alpha}=\nrm{Q_j\gamma}\le|Q_j\gamma-P_j|=|\delta_j|=q^{-s}.
\]
So $\deg g<L:=n+\kappa-s+1$ and $\nrm{g\alpha}<q^{-M}$ with $M:=s-1$. Here $L+M=n+\kappa$ and $L-M=n+\kappa-2s+2$. Given $g$, the polynomial $B$ is a monic divisor of $g$, so the number of $B$ is at most $q^{o(n)}\mathcal N(L,M)$. By Lemma~\ref{lem:count}, this is at most $q^{o(n)}q^{L-M+E/2}$ if $L\ge M$, and at most $q^{o(n)}q^{(E+L-M)/2}$ if $L<M$.

Since $s>\frac{k}{1+\varepsilon}\ge k(1-\varepsilon)$, we have $2s\ge2k-2\varepsilon k$ and $k-s<\varepsilon k$. Multiplying by $|\F_B|=q^{k-\kappa}$, we find that the contribution of these $B$ is at most
\[
q^{o(n)}\Bigl(q^{k-\kappa+n+\kappa-2s+2+E/2}+q^{k-\kappa+(E+n+\kappa-2s+2)/2}\Bigr)
\le q^{o(n)}\Bigl(q^{n-k+2\varepsilon k+E/2}+q^{(n-\kappa+E)/2+\varepsilon k}\Bigr).
\]
As $k>n/2$, this is bounded by the last two terms of the claim.

\emph{(c) $j\ge0$ and $s<\kappa$.} Then $s<k$, and $\F_B=Q_{j+1}\M_{k-s}$ by Lemma~\ref{lem:fibre}. Since $\deg Q_{j+1}=s<k\le n$, Lemma~\ref{lem:elem}(f) gives $|\sum_{D\in\F_B}\mu(D)|\le q^{1+o(n)}$. Summing over the $q^{n-k}$ values of $B$ gives at most $q^{n-k+1+o(n)}\le q^{n/2+o(n)}$, which is bounded by the second term of the claim.
\end{proof}

\subsection{Conclusion}
By~\eqref{eq:rangeE}, the exponents occurring in the error terms satisfy
\begin{align*}
\tfrac{n+E}2&\le n-\kappa-3\varepsilon n,\\
\tfrac n2+2\varepsilon n+\tfrac E2&\le n-\kappa-\varepsilon n,\\
\tfrac{n-\kappa+E}2+\varepsilon n&\le n-\kappa-2\varepsilon n.
\end{align*}
The first two inequalities are equivalent to~\eqref{eq:rangeE}. For the third, we use $E\le n-2\kappa-6\varepsilon n\le n-\kappa-6\varepsilon n$, which gives $\frac{n-\kappa+E}2\le n-\kappa-3\varepsilon n$.

We now insert~\eqref{eq:small}, with $K=\lfloor n/2\rfloor$, $E'=E$ and $q^{K+E/2}\le q^{(n+E)/2}$, and Proposition~\ref{prop:large}, for $\lfloor n/2\rfloor<k\le n$ with weight $k\le n$, into~\eqref{eq:decomp}. We get
\[
S=q^{n-\kappa}+O\bigl(n^2q^{n-\kappa-\min(\psi\varepsilon/(3p),\,\varepsilon)n+o(n)}\bigr).
\]
Since $\psi<\frac12$ and $p\ge3$, we have $\psi\varepsilon/(3p)<\varepsilon$. This proves Theorem~\ref{thm:main} with any $c<\psi\varepsilon/(3p)$. \qed

\section{Proof of Theorem~\ref{thm:main2}}\label{sec:proof2}
The proof follows the same pattern as in Section~\ref{sec:proof}. The difference is that the large divisors are now treated entirely by cancellation, using both directions of the fibres (Proposition~\ref{prop:twodir}). As a result, $E_\alpha(n+\kappa)$ enters only through the counting in the small-divisor range, and the large-divisor range imposes no condition on $\alpha$ at all. In Section~\ref{sec:proof}, by contrast, the large-divisor range is only available for $k>n/2$. The small-divisor estimate~\eqref{eq:small} therefore has to be used up to $K=\lfloor n/2\rfloor$, where it costs $q^{(n+E)/2}$, and case~(b) of Proposition~\ref{prop:large} gives a term of the same size. Here the large-divisor range starts at $k\approx\kappa/(1-2\eta)$, which allows $K\approx n-\kappa-\frac E2$. This is why~\eqref{eq:range2} contains $2\kappa+\frac12E$ in place of the quantity $2\kappa+E$ in~\eqref{eq:range}. The price is the factor~$4$ in~\eqref{eq:qcond2}.

Fix $\alpha$, $n$ and $\kappa$ with $\kappa\ge n/3$ and~\eqref{eq:range2}, and put $E=E_\alpha(n+\kappa)$ and $\V=\V_n(\alpha,\kappa)$. Since $E\ge1$, condition~\eqref{eq:range2} implies $\kappa\le\frac n2$. Together with $\kappa\ge n/3$, it also gives
\begin{equation}\label{eq:E2}
E\le2(1-4\eta)n-4\kappa\le n-\kappa-8\eta n,
\end{equation}
where the second inequality is equivalent to $n\le3\kappa$.

\emph{Size of $\V$.} We apply Lemma~\ref{lem:count}(ii) with $L=\kappa$ and $M=n$, so that $L+M=n+\kappa$. Since $M-L=n-\kappa\ge E$ by~\eqref{eq:E2}, we obtain $|\V|=q^{n-\kappa}$ exactly as in \S\ref{sec:sizeV}. We use the decomposition~\eqref{eq:decomp}.

\emph{Small divisors.} Let $K=\lfloor n-\kappa-\frac E2-\eta n\rfloor$, so that $1\le K\le n-\kappa$; the lower bound is checked below. Since only the scale $L+M=n+\kappa$ occurs in Proposition~\ref{prop:small}, we may use~\eqref{eq:small} with $E'=E=E_\alpha(n+\kappa)$. The two exponents in its error term satisfy $K+\frac E2\le n-\kappa-\eta n$ and, by~\eqref{eq:E2}, $\frac{n-\kappa+E}2\le n-\kappa-4\eta n$. Hence
\[
-\sum_{k\le K}k\,\Sigma_k=q^{n-\kappa}+O\bigl(q^{n-\kappa-\eta n+o(n)}\bigr).
\]

\emph{Large divisors.} By~\eqref{eq:range2} we have $n-\kappa-\frac E2\ge\kappa+4\eta n$, and therefore $K\ge\kappa+3\eta n-1$. Since $\kappa\le\frac n2$ and $\eta\le\frac14$, we also have $\frac{\kappa}{1-2\eta}\le\kappa(1+4\eta)\le\kappa+2\eta n$. Hence, for $n$ large, every $k>K$ satisfies
\begin{equation}\label{eq:klarge}
k\ge\frac{\kappa}{1-2\eta}\qquad\text{and}\qquad k\ge\kappa+\eta n .
\end{equation}

\begin{proposition}\label{prop:large2}
For every $k\le n$ satisfying~\eqref{eq:klarge},
\[
|\Sigma_k|\ll q^{\,n-\kappa-c_0\eta^2n}+q^{\,n-\kappa-\eta n+o(n)},
\]
with $c_0=c_0(q,\eta)$ from Proposition~\ref{prop:twodir} and an implied constant depending only on $q$ and $\eta$. No condition on $\alpha$ is used.
\end{proposition}

\begin{proof}
As before, $\Sigma_k=\sum_{B\in\M_{n-k}}\sum_{D\in\F_B}\mu(D)$. For each $B$ we let $\gamma=B\alpha$, $j=j_\gamma(k+\kappa)$ and $s=d_{j+1}(\gamma)$, and we distinguish three cases.

If $j=-1$, then $\F_B=\M_k$, and the inner sum vanishes since $k\ge2$. If $j\ge0$ and $s>k$, then $\F_B=\emptyset$.

If $j\ge0$ and $s<\kappa$, then $\F_B=Q_{j+1}\M_{k-s}$, and the inner sum is at most $q^{1+o(n)}$ by Lemma~\ref{lem:elem}(f). By~\eqref{eq:klarge}, these $B$ contribute at most $q^{n-k+1+o(n)}\le q^{n-\kappa-\eta n+o(n)}$.

Finally, let $j\ge0$ and $\kappa\le s\le k$. By the first condition in~\eqref{eq:klarge}, $(k-s)+(s-\kappa)=k-\kappa\ge2\eta k$, so that $\max(k-s,s-\kappa)\ge\eta k$. Proposition~\ref{prop:twodir} therefore bounds the inner sum by $Cq^{k-\kappa-c_0\eta k}$, and $k\ge\eta n$. Summing over the $q^{n-k}$ values of $B$ gives $Cq^{n-\kappa-c_0\eta^2n}$.
\end{proof}

Summing Proposition~\ref{prop:large2} over $K<k\le n$ with weights $k\le n$, and adding the contribution of the small divisors, we obtain $S=q^{n-\kappa}+O(n^2q^{n-\kappa-cn})$ for any $c<\min(c_0\eta^2,\eta)$. This proves Theorem~\ref{thm:main2}. \qed

\section{Proofs of the corollaries}\label{sec:cor}
The corollaries are obtained by choosing $\kappa=\lceil\theta n\rceil$ and checking that the hypotheses of the theorems are satisfied for the relevant $n$. For Corollary~\ref{cor:all} we balance Theorem~\ref{thm:weil}(b), which is strong when $\alpha$ has large partial quotients, against Corollary~\ref{cor:sigma}, which is strong when it has small ones.

\begin{proof}[Proof of Corollary~\ref{cor:tau}]
(a) Since $\tau(\alpha)<\tau$, there is $C_\alpha$ with $\deg a_{i+1}\le\tau d_i+C_\alpha$ for all $i$, and hence $E_\alpha(X)\le A_\alpha(X)\le\tau X+C_\alpha$. With $\kappa=\lceil\theta n\rceil$, the left side of~\eqref{eq:range} is therefore at most
\[
\theta n+1+\tfrac12\bigl(\tau(n+\theta n+1)+C_\alpha\bigr)=\tfrac12\bigl(\theta(2+\tau)+\tau\bigr)n+O(1).
\]
Since $\theta(2+\tau)+\tau<1-6\varepsilon$, this is at most $(\frac12-3\varepsilon)n$ for large $n$.

Theorem~\ref{thm:main} and Lemma~\ref{lem:elem}(e) now give~\eqref{eq:primecount}, since $n/2<n-\kappa-cn$. Every $P$ counted satisfies $\nrm{\alpha P}<q^{-\kappa}\le|P|^{-\theta}$, and so $\theta(\alpha)\ge\theta$. The last statement follows by letting $\tau\to0$: if $\tau(\alpha)=0$, then for every $\theta<\frac12-3\varepsilon$ there is $\tau>0$ with $\theta(2+\tau)<1-\tau-6\varepsilon$.

(b) Note that $\theta<\frac12$. Choose $\sigma'\in(\sigma(\alpha),\sigma)$. Exactly as in the proof of Corollary~\ref{cor:sigma} below, there are infinitely many $j$ and $n$ such that
\[
X:=n+\lceil\theta n\rceil\in[d_j+d_{j+1},\,d_{j+1}+d_{j+2}),\qquad j_\alpha(X)=j,\qquad E_\alpha(X)\le\sigma'X .
\]
With $\kappa=\lceil\theta n\rceil$, the left side of~\eqref{eq:range} is at most
\[
\theta n+1+\tfrac{\sigma'}2(n+\theta n+1)=\tfrac12\bigl(\theta(2+\sigma')+\sigma'\bigr)n+O(1).
\]
Since $\theta(2+\sigma')+\sigma'<\theta(2+\sigma)+\sigma<1-6\varepsilon$, this is at most $(\frac12-3\varepsilon)n$ for large $j$. Theorem~\ref{thm:main} gives~\eqref{eq:primecount} for these $n$, and we conclude as in~(a).
\end{proof}

\begin{proof}[Proof of Corollary~\ref{cor:sigma}]
Choose $\sigma'\in(\sigma(\alpha),\sigma)$. By the definition of $\sigma(\alpha)$, there are infinitely many $j$ with $\max(\deg a_{j+1},\deg a_{j+2})\le\sigma'(d_j+d_{j+1})$. Fix such a large $j$. The interval $I_j=[d_j+d_{j+1},\,d_{j+1}+d_{j+2})$ contains $d_{j+2}-d_j\ge2$ consecutive integers. The map $n\mapsto n+\lceil\theta n\rceil$ is increasing and unbounded, and it increases in steps of $1$ or $2$ because $\theta<1$. So it cannot jump over $I_j$, and there is $n$ with $X:=n+\lceil\theta n\rceil\in I_j$. Then $j_\alpha(X)=j$, and
\[
E_\alpha(X)=\max(\deg a_{j+1},\deg a_{j+2})\le\sigma'(d_j+d_{j+1})\le\sigma'X.
\]

With $\kappa=\lceil\theta n\rceil\ge n/3$, the left side of~\eqref{eq:range2} is at most
\[
2\theta n+2+\tfrac{\sigma'}2(n+\theta n+1)=\tfrac12\bigl(\theta(4+\sigma')+\sigma'\bigr)n+O(1).
\]
Since $\sigma'<\sigma$, we have $\theta(4+\sigma')+\sigma'<\theta(4+\sigma)+\sigma<2-8\eta$, so this is at most $(1-4\eta)n$ for large $j$. Theorem~\ref{thm:main2} gives~\eqref{eq:primecount} for this $n$. As $j\to\infty$, we have $X\ge d_j+d_{j+1}\to\infty$, so we obtain infinitely many $n$. The last statement follows by letting $\sigma\to0$.
\end{proof}

\begin{proof}[Proof of Corollary~\ref{cor:all}]
Let $\frac13\le\theta<\sqrt2-1$, and put
\[
f(\tau)=\frac{1+\tau}{3+\tau},\qquad\sigma(\tau)=\frac{\tau(1+\tau)}{2+\tau},\qquad h(\tau)=\frac{2-\sigma(\tau)}{4+\sigma(\tau)} .
\]
Thus $f(\tau)$ is the exponent of Theorem~\ref{thm:weil}(b), and $h(\tau)$ that of Corollary~\ref{cor:sigma} with $\sigma=\sigma(\tau)$. The function $f$ is increasing on $[0,1]$, and $\sigma$ is increasing, so $h$ is decreasing. At $\tau=\sqrt2-1$ the two exponents agree: $f=\frac{\sqrt2}{2+\sqrt2}=\sqrt2-1$; moreover $\sigma=\frac{(\sqrt2-1)\sqrt2}{\sqrt2+1}=3\sqrt2-4$, so $h=\frac{6-3\sqrt2}{3\sqrt2}=\sqrt2-1$. Hence there is $\tau_1\in(0,1)$ with $f(\tau_1)>\theta$ and $h(\tau_1)>\theta$. The first condition is equivalent to $\tau_1>\frac{3\theta-1}{1-\theta}$, and the second to $2-\sigma_1-\theta(4+\sigma_1)>0$, so that $\eta$ as in the statement exists.

\emph{Case $\tau(\alpha)\ge\tau_1$.} Then $\min(\tau(\alpha),1)\ge\tau_1$, and Theorem~\ref{thm:weil}(b) gives $\theta(\alpha)\ge f(\tau_1)>\theta$.

\emph{Case $\tau(\alpha)<\tau_1$.} Then $\deg a_{i+1}\le\tau_1d_i$ for all large $i$. For large $j$ this gives $\deg a_{j+1}\le\tau_1d_j\le\tau_1d_{j+1}$ and $\deg a_{j+2}\le\tau_1d_{j+1}$, as well as $d_{j+1}=d_j+\deg a_{j+1}\le(1+\tau_1)d_j$. Hence
\[
\frac{\max(\deg a_{j+1},\deg a_{j+2})}{d_j+d_{j+1}}\le\frac{\tau_1d_{j+1}}{d_{j+1}\bigl(\frac1{1+\tau_1}+1\bigr)}=\frac{\tau_1(1+\tau_1)}{2+\tau_1}=\sigma_1 ,
\]
and so $\sigma(\alpha)\le\sigma_1$. Choose $\sigma>\sigma_1$ so close to $\sigma_1$ that $\theta(4+\sigma)<2-\sigma-8\eta$. Then Corollary~\ref{cor:sigma} gives $\theta(\alpha)\ge\theta$.

\emph{The alternative bound.} Let $\frac13\le\theta<\frac{\sqrt{17}-1}8$, and put $g(\tau)=\frac{1-\sigma(\tau)}{2+\sigma(\tau)}$, the exponent of Corollary~\ref{cor:tau}(b) with $\sigma=\sigma(\tau)$. Since $1-\sigma(\tau)=\frac{2-\tau^2}{2+\tau}$ and $2+\sigma(\tau)=\frac{\tau^2+3\tau+4}{2+\tau}$, we have $g(\tau)=\frac{2-\tau^2}{\tau^2+3\tau+4}$, which is decreasing on $[0,1]$. The equation $f(\tau)=g(\tau)$ reads
\[
2\tau^3+7\tau^2+5\tau-2=(\tau+2)(2\tau^2+3\tau-1)=0,
\]
so $f$ and $g$ meet at $\tau^*=\frac{\sqrt{17}-3}4=0.2807\ldots$, with common value $f(\tau^*)=\frac{\sqrt{17}-1}8$. Hence there is $\tau_1\in(0,1)$ with $f(\tau_1)>\theta$ and $g(\tau_1)>\theta$. The second condition is $(1-\sigma_1)-\theta(2+\sigma_1)>0$, so that $\varepsilon$ as in the statement exists; moreover $\varepsilon<\frac1{18}\le\frac1{10}$, because $\theta\ge\frac13$.

The case $\tau(\alpha)\ge\tau_1$ is handled by Theorem~\ref{thm:weil}(b) as before. If $\tau(\alpha)<\tau_1$, then $\sigma(\alpha)\le\sigma_1$, as shown above. Choose $\sigma>\sigma_1$ so close to $\sigma_1$ that $\theta(2+\sigma)<1-\sigma-6\varepsilon$; then Corollary~\ref{cor:tau}(b) gives $\theta(\alpha)\ge\theta$.
\end{proof}

\section{Limits of the method}\label{sec:remarks}

In this final section we \chg{discuss} where the method stops, and what seems to be needed to go further.

\subsection{\texorpdfstring{Why the exponent $\frac12$}{Why the exponent 1/2}}
The small-divisor range is exact as long as $\chg{k\le n-\kappa}$, and the large-divisor range needs fibres $\F_B$ that contain progressions of positive relative length, which requires $k>\kappa$. Both are needed at $k=n/2$, which forces $\kappa<n/2$. For $\kappa>n/2$ and $n-\kappa<k<\kappa$, Lemma~\ref{lem:fibre} shows that every nonempty fibre is $\M_k$ or $Q_{j+1}\M_{k-s}$ with $s\le k$. The fibres with $s\le k-\eta n$ have negligible M\"obius sums by Lemma~\ref{lem:elem}(f). Those with $s$ close to $k$ consist of a few multiples of a single convergent denominator $Q_{j+1}(B\alpha)$, and their contribution is essentially
\[
\sum_{B}\mu\bigl(Q_{j+1}(B\alpha)\bigr),
\]
a sum of $\mu$ over convergent denominators of the multiples $B\alpha$. We know no estimate for such sums.

\subsection{A Fourier-dual formulation and short trace sums}
The following heuristic discussion explains why a Fourier-dual approach does not, by itself, go beyond $\frac12$; we do not make it rigorous. One can also write $\V_n$ as a union of residue classes modulo a convergent denominator $f$ of $\alpha$ and detect the classes by additive characters. For $k$ in the range above, $\Sigma_k$ then becomes a weighted sum, over about $q^{\kappa}$ frequencies, of sums of $\mu$ over Diophantine sets $\V_k(h/f,n-k)$, twisted by additive phases. These sets split into long progressions as in Lemma~\ref{lem:fibre}, and the phases are trace functions with slopes at most $1$ at infinity. So one may hope to apply the bounds of Sawin and Shusterman for short sums of such trace functions~\cite{SS3}, and~\cite[Theorems~1.10 and~5.1]{SS2}, to the individual sums. However, a count of sizes suggests that, after normalisation, the total size of these sums is about $q^{k}$, while the target is $q^{n-\kappa}$. Even square-root cancellation along the progressions, a saving of about $q^{(k-n/2)/2}$, only compensates for this when $k<\frac32n-2\kappa$, and at $k=\kappa$ this again requires $\kappa<n/2$.

\subsection{\texorpdfstring{Squarefree moduli and the size of $q$}{Squarefree moduli and the size of q}}
Our moduli $Q_j(B\alpha)$ and $Q_{j+1}(B\alpha)$ are arbitrary. The inputs that allow this,~\cite[Theorem~4.5]{SS} and~\cite[Theorem~5.1]{SS2}, rely on the $p$-th power structure of $\mu$ in characteristic $p$, which is the source of the condition $q\gg p^2$; in particular, $q$ prime is never allowed. For squarefree moduli, Sawin~\cite[Theorems~1.1, 1.2 and 1.7]{Sawin} proves level of distribution results for $\mu$, $\Lambda$ and all factorization functions, with conditions on $q$ independent of $p$ and levels tending to $1$. Two consequences:
\begin{enumerate}[label=(\alph*)]
\item If one could choose squarefree moduli in Lemma~\ref{lem:fibre} for all but an exponentially small proportion of $B$, one could hope for a version of Theorem~\ref{thm:main} with a condition on $q$ independent of $p$, including prime $q$. This would also need a version of~\cite[Theorem~1.7]{Sawin} for progressions intersected with short intervals. We do not know how to do this.
\item If $\alpha$ has infinitely many squarefree convergent denominators (or squarefree members of the pencils $q_m+cq_{m-1}$), Theorem~\ref{thm:sqfree} uses~\cite[Theorem~1.2]{Sawin} directly and gives $\theta(\alpha)\ge2\lambda(q)-1$, which exceeds $\frac12$ for $q>23172$, for every $q$ including prime $q$.
\end{enumerate}

\section*{Acknowledgements}
\chg{The author is grateful to Prof.~Stephan Baier for suggesting the approach based on averages of primes over segments of arithmetic progressions and for many helpful discussions. The author also thanks Prof.~Will Sawin for helpful email correspondence.} This work was started at the Indian Statistical Institute, Kolkata, and completed at the Indian Institute of Technology Bombay; the author thanks both institutes for excellent working environments. The author was supported by a Research Associate fellowship at the Indian Statistical Institute, Kolkata, and is currently supported as a Project Research Scientist at the Indian Institute of Technology Bombay (Project Code No.~RD/0525-IRCCSH0-005).

\end{document}